\documentclass[10pt]{article}
\usepackage [dvips]{graphicx}
\usepackage{amssymb}
\newcommand{\devp}[2]{ \frac{\partial #1}{\partial #2}}

\newcommand{\der}[2]{ \frac{d #1}{d #2}}

\newcommand{\sv}{\hspace{1pt} ,}

\newcommand{\eq}[1]{~\mbox{{\rm(\ref{#1})}}}

\newcommand{\arrayvspace}{\\[0.2in]}
\newcommand{\be}{\begin{equation}}
\newcommand{\ee}{\end{equation}}
\newcommand{\bea}{\begin{eqnarray}}
\newcommand{\eea}{\end{eqnarray}}
\newcommand{\bean}{\begin{eqnarray*}}
\newcommand{\eean}{\end{eqnarray*}}
\def\ba{\begin{array}{l}\displaystyle}
\def\ea{\end{array}}

\def\F{{\cal F}}
\def\r{{r}}
\def\h{{\Delta t}}
\def\hh{{\Delta x}}

\newcommand{\rf}[1]{(\ref{#1})}

\newcommand{\tA}{\tilde{A}}
\newcommand{\ta}{\tilde{a}}
\newcommand{\tw}{\tilde{w}}
\newcommand{\tc}{\tilde{c}}

\newcommand{\Rz}{{\cal R}}
\newcommand{\e}{{\varepsilon}}
\newcommand{\R}{{\mathbb R}}
\newcommand{\E}{{\cal E}}

\newtheorem{definition}{Definition}[section]
\newtheorem{lemma}{Lemma}[section]
\newtheorem{theorem}{Theorem}[section]

\def\proof{{\bf Proof}:\\}
\def\endproof{{\par \raggedleft \vrule height4pt width 4pt \par}}

\newenvironment{remarks}{{\flushleft \bf Remarks:}}{}

\begin{document}

\title{Implicit-explicit Runge-Kutta schemes and applications to hyperbolic systems with relaxation\thanks{This work was
supported by the European network HYKE, funded by the EC as
contract HPRN-CT-2002-00282.}}

% The thanks line in the title should be filled in if there is
% any support acknowledgement for the overall work to be included
% This \thanks is also used for the received by date info, but
% authors are not expected to provide this.

\author{Lorenzo Pareschi\thanks{Department of Mathematics, University of Ferrara, Via
Machiavelli 35, I-44100 Ferrara, Italy. ({\tt
pareschi@dm.unife.it})}\and Giovanni Russo\thanks{Department of
Mathematics and Computer Science, University of Catania, Via
A.Doria 6, 95125 Catania, Italy. ({\tt russo@dmi.unict.it})}}

\date{May 6, 2004}

\maketitle

\begin{abstract}
We consider new implicit-explicit (IMEX) Runge-Kutta methods for
hyperbolic systems of conservation laws with stiff relaxation
terms. The explicit part is treated by a
strong-stability-preserving (SSP) scheme, and the implicit part is
treated by an L-stable diagonally implicit Runge-Kutta methods
(DIRK). The schemes proposed are asymptotic preserving (AP) in the
zero relaxation limit. High accuracy in space is obtained by
Weighted Essentially Non Oscillatory (WENO) reconstruction. After
a description of the mathematical properties of the schemes,
several applications will be presented.
\end{abstract}

{\bf Keywords :} Runge-Kutta methods, hyperbolic systems with
relaxation, stiff systems, high order shock capturing schemes.

\medskip
{\bf AMS Subject Classification :} 65C20, 82D25

%\tableofcontents

\section{Introduction}          \label{sc:introduction}

Several physical phenomena of great importance for applications
are described by stiff systems of differential equations in the
form \be
  \partial_t U  =  \F(U) + \frac1{\varepsilon}R(U),
  \label{eq:stiff-ODE1}
\ee where $U=U(t)\in \R^N$, $\F, R:\R^N\to \R^N$ and $\varepsilon >
0$ is the stiffness parameter.

System \rf{eq:stiff-ODE1} may represent a system of $N$ ODE's or a
discretization of a system of PDE's, such as, for example,
convection-diffusion equations, reaction-diffusion equations and
hyperbolic systems with relaxation.

In this work we consider the latter case, which in recent years
has been a very active field of research, due to its great impact
on applied sciences \cite{CLL, Liu}. For example, we mention that
hyperbolic systems with relaxation appears in kinetic theory of
rarefied gases, hydrodynamical models for semiconductors,
viscoelasticity, multiphase flows and phase transitions, radiation
hydrodynamics, traffic flows, shallow waters, etc.

In one space dimension these systems have the form
\be
\partial_t U + \partial_x F(U) = \frac1{\e} R(U),\quad x \in \R,
\label{eq:hsr}
\ee
that corresponds to (\ref{eq:stiff-ODE1}) for
$\F(U)=-\partial_x F(U)$. In (\ref{eq:hsr}) the Jacobian matrix
$F'(U)$ has real eigenvalues and admits a basis of eigenvectors
$\forall\,U\in \R^N$ and $\e$ is called {\em relaxation
parameter}.

The development of efficient numerical schemes for such systems is
challenging, since in many applications the relaxation time varies
from values of order one to very small values if compared to the
time scale determined by the characteristic speeds of the system.
In this second case the hyperbolic system with relaxation is said
to be stiff and typically its solutions are well approximated by
solutions of a suitable reduced set of conservation laws called
{\em equilibrium system} \cite{CLL}.

Usually it is extremely difficult, if not impossible, to split the
problem in separate regimes and to use different solvers in the
stiff and non stiff regions. Thus one has to use the original
relaxation system in the whole computational domain. The
construction of schemes that work for all ranges of the relaxation
time, using coarse grids that do not resolve the small relaxation
time, has been studied mainly in the context of upwind methods
using a method of lines approach combined with suitable operator
splitting techniques \cite{CJR, Jin-RK} and more recently in the
context of central schemes \cite{LRR, Pa}.

Splitting methods have been widely used for such problems. They
are attractive because of their simplicity and robustness. Strang
splitting provides second order accuracy if each step is at least
second order accurate \cite{Strang}. This property is maintained
under fairly mild assumptions even for stiff problems
\cite{Lubich}. However, Strang splitting applied to hyperbolic
systems with relaxation reduces to first order accuracy when the
problem becomes stiff. The reason is that the kernel of the
relaxation operator is non trivial, which corresponds to a
singular matrix in the linear case, and therefore the assumptions
in \cite{Lubich} are not satisfied.

Furthermore with a splitting strategy it is difficult to obtain
higher order accuracy even in non stiff regimes (high order
splitting schemes can be constructed, see \cite{Dia-Shatzman2},
but they are seldom used because of stability problems).

Recently developed Runge-Kutta schemes overcome this difficulties,
providing basically the same advantages of the splitting schemes,
without the drawback of the order restriction
\cite{CJR,Jin-RK,Zhong}.

In this paper we will present some recent results on the
development of high order, underresolved Runge-Kutta time
discretization for such systems. In particular, using the
formalism of implicit-explicit (IMEX) Runge-Kutta schemes
\cite{ARS, ARW, PR, CK, Zhong} we derive new IMEX schemes up to order 3
that are strong-stability-preserving (SSP) for the limiting system
of conservation laws (such methods were originally referred to as
total variation diminishing (TVD) methods, see \cite{GS, GST,
SR}).

To this aim, we derive general conditions that guarantee the
asymptotic preserving property, i.e.\ the consistency of the
scheme with the equilibrium system and the asymptotic accuracy,
i.e. the order of accuracy is maintained in the stiff limit. For a
stability analysis of IMEX schemes we refer to \cite{PR3}.

The rest of the paper is organized as follows. In Section 2 we
introduce the general structure of IMEX Runge-Kutta schemes.
Section 3 is devoted to IMEX Runge-Kutta schemes applied to
hyperbolic systems with relaxation. In Section 4 we describe space
discretization obtained by conservative finite-volume and finite
difference schemes. In Section 5 we present some numerical results
concerning the accuracy of IMEX schemes when applied to a
prototype hyperbolic system with relaxation. Finally in Section 6
we investigate the performance of the schemes in several realistic
applications to shallow waters, traffic flows and granular gases.

%------------------------------------------------------------------------------

\section{IMEX Runge-Kutta schemes}
\label{sec:2}

 An
IMEX Runge-Kutta scheme consists of applying an implicit
discretization to the source terms and an explicit one to the
non stiff term. When applied to system (\ref{eq:stiff-ODE1}) it takes the form
\begin{eqnarray}
 U^{(i)} & = & U^n - \h \sum_{j=1}^{i-1} \ta_{ij} \partial_x F(U^{(j)}) +
                 \h \sum_{j=1}^\nu a_{ij} \frac{1}{\varepsilon} R(U^{(j)}),
                                                                 \label{eq:RKEI1}\\
 U^{n+1} & = & U^n - \h \sum_{i=1}^{\nu} \tw_{i} \partial_x F(U^{(i)})  +
                 \h \sum_{i=1}^\nu w_{i} \frac{1}{\varepsilon} R(U^{(i)}).
                                                                 \label{eq:RKEI2}
\end{eqnarray}
The matrices $\tA = (\ta_{ij})$, $\ta_{ij}=0$ for $j \geq i$ and
$A = (a_{ij})$ are $\nu \times \nu$ matrices such that the
resulting scheme is explicit in $F$, and implicit in $R$. An IMEX
Runge-Kutta scheme is characterized by these two matrices and the
coefficient vectors $\tw= (\tw_1,\ldots,\tw_\nu)^T$, $w =
(w_1,\ldots,w_\nu)^T$.

Since the simplicity and efficiency of solving the algebraic
equations corresponding to the implicit part of the discretization
at each step is of paramount importance, it is natural to consider
diagonally implicit Runge-Kutta (DIRK) schemes \cite{HW} for the
source terms ($a_{ij}=0$, for $j > i$).

IMEX Runge-Kutta schemes can be represented by a double {\em tableau\/} in the usual
Butcher notation,
\[
\begin{array}{c|c}
              \tc & \tA \\
              \hline \\
              & \tw^T
\end{array}
\qquad\qquad
\begin{array}{c|c}
              c & A \\
              \hline \\
              & w^T
\end{array}
\]
where the coefficients $\tc$ and $c$ used for the treatment of non
autonomous systems, are given by the usual relation \be
   \tc_i = \sum_{j=1}^{i-1} \ta_{ij}, \quad
     c_i = \sum_{j=1}^{i}   a_{ij}.
\label{eq:C} \ee
The use of a DIRK scheme for $R$ is a sufficient condition to
guarantee that $F$ is always evaluated explicitly.

In the case of system (\ref{eq:hsr}), in order to obtain a
numerical scheme, a suitable space discretization of equations
(\ref{eq:RKEI1})-(\ref{eq:RKEI2}) is required. This discretization
can be performed at the level of the original system
(\ref{eq:hsr}) in which case one has a system of ODEs that is then
solved in time by the IMEX scheme (method of lines). Alternatively
one can apply a suitable space discretization directly to the time
discrete equations (\ref{eq:RKEI1})-(\ref{eq:RKEI2}).

Finally we remark that previously developed schemes, such as the
Additive semi-implicit Runge-Kutta methods of Zhong \cite{Zhong},
and the splitting Runge-Kutta methods of Jin et al. \cite{Jin-RK},
\cite{CJR} can be rewritten as IMEX-RK schemes \cite{PR3}.

\subsection{Order conditions}

The general technique to derive order conditions for Runge-Kutta
schemes is based on the Taylor expansion of the exact and
numerical solution.

In particular, conditions for schemes of order $p$ are obtained by imposing
that the solution of system (\ref{eq:hsr}) at time $t=t_0+\Delta t$, with a
given initial condition at time $t_0$,
 agrees with the numerical
solution obtained by one step of a Runge-Kutta scheme with the
same initial condition, up to order $\h^{p+1}$.

Here we report the order conditions for IMEX Runge-Kutta schemes
up to order $p=3$, which is already considered high order for PDEs
problems.

We apply scheme (\ref{eq:RKEI1})-(\ref{eq:RKEI2}) to system
(\ref{eq:hsr}), with $\varepsilon=1$. We assume that the
coefficients $\tc_i$, $c_i$, $\ta_{ij}$, $a_{ij}$ satisfy
conditions (\ref{eq:C}). Then the order conditions are the
following
\paragraph*{First order}
\be
   \sum_{i=1}^\nu \tw_i = 1, \quad \sum_{i=1}^\nu w_i = 1.
   \label{eq:cond1}
\ee
\paragraph*{Second order}
\be
   \sum_i \tw_i \tc_i = 1/2, \quad \sum_i w_i c_i = 1/2,
   \label{eq:cond2}
\ee
\be
   \sum_i \tw_i c_i = 1/2, \quad \sum_i w_i \tc_i = 1/2.
   \label{eq:mixed2}
\ee
\paragraph*{Third order}
\be
   \sum_{ij} \tw_i \ta_{ij} \tc_j = 1/6, \quad
   \sum_{i}  \tw_i \tc_{i}  \tc_i = 1/3, \quad
   \sum_{ij}   w_i   a_{ij}   c_j = 1/6, \quad
   \sum_{i}    w_i   c_{i}    c_i = 1/3,
\label{eq:cond3}
  \ee
\bea
   \begin{array}{c}
   \sum_{ij} \tw_i \ta_{ij}   c_j = 1/6, \quad
   \sum_{ij} \tw_i   a_{ij} \tc_j = 1/6, \quad
   \sum_{ij} \tw_i   a_{ij}   c_j = 1/6, \arrayvspace
   \sum_{ij}   w_i \ta_{ij}   c_j = 1/6, \quad
   \sum_{ij}   w_i   a_{ij} \tc_j = 1/6, \quad
   \sum_{ij}   w_i \ta_{ij} \tc_j = 1/6, \arrayvspace
   \sum_{i}  \tw_i   c_{i}    c_i = 1/3, \quad
   \sum_{i}  \tw_i \tc_{i}    c_i = 1/3, \quad
   \sum_{i}    w_i \tc_{i}  \tc_i = 1/3, \quad
   \sum_{i}    w_i \tc_{i}    c_i = 1/3.
   \end{array}                                     \label{eq:mixed3}
\eea
Conditions \rf{eq:cond1}, \rf{eq:cond2}, \rf{eq:cond3} are
the standard order conditions for the two {\em tableau\/}, each of
them taken separately. Conditions \rf{eq:mixed2} and
\rf{eq:mixed3} are new conditions that arise because of the
coupling of the two schemes.

The order conditions will simplify a lot if $\tilde{c} = c$. For
this reason only such schemes are considered in \cite{ARS}. In
particular, we observe that, if the two tableau differ only for
the value of the matrices $A$, i.e.\ if $\tc_i = c_i$ and $\tw_i =
w_i$, then the standard order conditions for the two schemes are
enough to ensure that the combined scheme is third order. Note,
however, that this is true only for schemes up to third order.

Higher order conditions can be derived as well using a
generalization of Butcher 1-trees to 2-trees, see \cite{CK}.
However the number of coupling conditions increase dramatically
with the order of the schemes. The relation between coupling conditions
and accuracy of the schemes is reported in Table \ref{tab:order}.

\begin{table}
\begin{center}
\begin{tabular}{c|c|c|c|c}
{\bf IMEX-RK} & \multicolumn{4}{c}{\bf Number of coupling conditions}\\
{\bf order}  & General case  & $\tw_i = w_i$ & $\tilde{c} = c$ &
$\tilde{c} = c$ and $\tw_i =
w_i$\\
\hline
1 & 0 & 0 & 0 & 0\\
2 & 2 & 0 & 0 & 0\\
3 & 12 & 3 & 2 & 0\\
4 & 56 & 21 & 12 & 2\\
5 & 252 & 110 & 54 & 15\\
6 & 1128 & 528 & 218 & 78
\end{tabular}
\end{center}
\caption{Number of coupling conditions in IMEX Runge-Kutta
schemes} \label{tab:order}
\end{table}

\section{Applications to hyperbolic systems with relaxation}

In this section we give sufficient conditions for asymptotic preserving and
asymptotic accuracy properties of IMEX schemes. This properties are strongly
related to L-stability of the implicit part of the scheme.

\subsection{Zero relaxation limit}
Let us consider here one-dimensional hyperbolic systems with
relaxation of the form (\ref{eq:hsr}). The operator
$R:\R^N\to\R^N$ is called a relaxation operator, and consequently
(\ref{eq:hsr}) defines a relaxation system, if there exists a
constant $n \times N$ matrix $Q$ with rank$(Q)=n<N$ such that \be
QR(U)=0\quad \forall\,\, U \in \R^N. \label{eq:cons} \ee This
gives $n$ independent conserved quantities $u=QU$. Moreover
%% Quello che hai scritto mi sembra una assunzione,
%% piu' che un fatto. Mi sembra piu' corretto scriveer come segue
we assume that equation $R(U)=0$ can be uniquely solved in terms
of $u$, i.e.
% such conserved quantities uniquely determine a local equilibrium value
\be U=\E(u) \,\,\, \hbox{such that} \,\,\, R(\E(u))=0.
\label{eq:maxwellian} \ee The image of $\E$ represents the
manifold of local equilibria of the relaxation operator $R$.

Using (\ref{eq:cons}) in (\ref{eq:hsr}) we obtain a system of $n$
conservation laws which is satisfied by every solution of
(\ref{eq:hsr})
\be
  \partial_t (Q U) + \partial_x (Q F(U)) = 0.
\label{eq:mom}
\ee

For vanishingly small values of the relaxation parameter $\e$ from
(\ref{eq:hsr}) we get $R(U)=0$ which by (\ref{eq:maxwellian})
implies $U=\E(u)$. In this case system (\ref{eq:hsr}) is well
approximated by the
%% Era stato chiamoato equilibrium affa fine di pag. 1.
equilibrium system \cite{CLL}.
\be
   \partial_t u + \partial_x G(u) = 0,
   \label{eq:fluid}
\ee
where $G(u)=QF(\E(u))$.

System (\ref{eq:fluid}) is the formal limit of system (\ref{eq:stiff-ODE1}) as
$\e \to 0$. The solution $u(x,t)$ of the system will be the limit of $QU$, with
$U$ solution of system (\ref{eq:stiff-ODE1}), provided a suitable condition on
the characteristic velocities of systems (\ref{eq:stiff-ODE1}) and
(\ref{eq:fluid}) is satisfied (the so called {\em subcharacteristic condition},
see \cite{Whitham, CLL}.)

\subsection{Asymptotic properties of IMEX schemes}

We start with the following
\begin{definition}
We say that an IMEX scheme for system \rf{eq:hsr} in the form
\rf{eq:RKEI1}-\rf{eq:RKEI2} is {\em asymptotic preserving\/} (AP)
if in the limit $\varepsilon\to 0$ the scheme becomes a consistent
discretization of the limit equation \rf{eq:fluid}.
\end{definition}
Note that this definition does not imply that the scheme preserves
the order of accuracy in $t$ in the stiff limit $\epsilon\to 0$.
In the latter case the scheme is said {\em asymptotically
accurate}.

%***** Add two examples of schemes that reduce or gain order, even if the order
%***** of explicit and implicit is the same.

%For example let us consider the results of a direct consistency
%analysis of the schemes of the previous section. Scheme SP(1,1,1),
%shown in Table \ref{tb:SP111}, is clearly AP. Scheme Jin(2,2,2) is
%AP, but it is not uniformly valid in $\varepsilon$, since the
%absolute-stability region of the implicit part of the scheme has a
%hole in the middle of the complex half plane $\Re z < 0$.
%% Mi pare che non riusciamo a mettere una referenza qui.
%% Se non sbaglio lo abbiamo trovato, ma non lo
%% abbiamo scritto ancora da nessuna parte
%% Per adesso non mettiamo alcuna referenza.
%Schemes Midpoint(1,2,2) and CN(2,2,2) are not AP even if both
%implicit parts of the schemes are A-stable \cite{CJR}. On the
%contrary, schemes CJR(3,2,2) and LRR(3,2,2) are AP and uniformly
%valid in $\varepsilon$, but only scheme LRR(3,2,2) preserves the
%second order accuracy in the stiff limit

In order to give sufficient conditions for the AP and
asymptotically accurate property, we make use of the following
simple
\begin{lemma}
If all diagonal elements of the triangular coefficient matrix $A$
that characterize the DIRK scheme are non zero, then
\be
  \lim_{\epsilon\to0} R(U^{(i)}) = 0. \label{eq:relax}
\ee
\end{lemma}

\proof In the limit $\epsilon\to 0$ from \rf{eq:RKEI1} we have
\[
  \sum_{j=1}^i a_{ij} R(U^j) = 0, \quad i=1,\ldots,\nu.
\]
Since the matrix $A$ is non-singular, this implies $R(U^i)=0,\>
i=1,\ldots,\nu$.
\endproof

  In order to apply the previous lemma, the vectors of $c$ and $\tilde{c}$
  cannot be equal. In fact $\tilde{c}_1 = 0$ whereas $c_1 \neq 0$.
  Note that if $c_1 = 0$ but $a_{ii} \neq 0$ for $i > 1$, then we still have $\lim_{\epsilon\to0} R(U^{(i)}) = 0$
  for $i>1$ but $\lim_{\epsilon\to0} R(U^{(1)}) \neq 0$ in general.
  The corresponding scheme may be inaccurate if the
  initial condition is not ``well prepared'' ($R(U_0) \neq 0$). In this case the scheme is not
  able to treat the so called ``initial layer'' problem and degradation of accuracy in
  the stiff limit is expected (see Section 5 and references \cite{CJR, PR,
  Pa}.) On the other hand, if the initial condition is ``well prepared'' ($R(U^{(0)})=0$),
  then relation \rf{eq:relax}, $i=1,\ldots,\nu$ holds even if $a_{11}=c_1=0$.

Next we can state the following
\begin{theorem}
  If $\det A \neq 0$ in the limit $\epsilon\to 0$, the IMEX
scheme \rf{eq:RKEI1}-\rf{eq:RKEI2} applied to system \rf{eq:hsr}
becomes the explicit RK scheme
  characterized by $(\tilde{A},\tilde{w},\tilde{c})$ applied to the limit
  equation \rf{eq:fluid}.
\label{th:main}
\end{theorem}

\proof
  As in the case of the continuous system, multiplying equations \rf{eq:RKEI1} and \rf{eq:RKEI2} by
  the matrix $Q$, and making use of the relation $QR(U) = 0 \> \forall\, U$, we obtain
\begin{equation}
  u^{(i)}  =  u_0 + h \sum_{j=1}^{i-1} \ta_{ij} QF(U^{(j)}),
\label{eq:RKEI1Q}
\end{equation}
and
\begin{equation}
  u_1 = u_0 + h \sum_{i=1}^{\nu} \tw_{i} QF(U^{(i)}),
\label{eq:RKEI2Q}
\end{equation}
where
\[
  u^{(i)} = QU^{(i)}, \> i=1,\ldots \nu, \quad u_0 = QU_0, \quad u_1 = QU_1.
\]
From the previous lemma, in the stiff limit one has
\[
  R(U^{(i)}) = 0, \> i=1,\ldots \nu,
\]
and, from property \rf{eq:maxwellian} of the relaxation operator
$R$, this implies
\[
  U^{(i)} = \E(u^{(i)}), \> i=1,\ldots \nu.
\]
Substituting this expression in \rf{eq:RKEI1Q}-\rf{eq:RKEI2Q} one
obtains
\begin{equation}
  u^{(i)}  =  u_0 + h \sum_{j=1}^{i-1} \ta_{ij} G(u^{(j)}),
\label{eq:RKEI1Q2}
\end{equation}
and
\begin{equation}
  u_1 = u_0 + h \sum_{i=1}^{\nu} \tw_{i} G(u^{(i)}),
\label{eq:RKEI2Q2}
\end{equation}
where $G(u) = QF(\E(u))$.
\endproof

Clearly one may claim that if the implicit part of the IMEX scheme
is A-stable or L-stable the previous theorem is satisfied. Note
however that this is true only if the {\em tableau} of the
implicit integrator does not contain any column of zeros that
makes it reducible to a simpler
%% Non scriverei ``but'', che in americano ha un
%% significato leggermente diverso.  Metterei nulla oppure
%% ``still'' o ``but still''
A-stable or L-stable form. Some remarks are in order.

\begin{remarks}

\begin{itemize}
\item[i)] There is a close analogy between hyperbolic systems with stiff relaxation
and differential algebraic equations (DAE) \cite{AP98}. The limit system
as $\epsilon\to 0$ is the analog of an index 1 DAE, in which the algebraic
equation is explicitly solved in terms of the differential variable.
In the context of DAE, the initial condition that we called ``well prepared'' is called ``consistent''.

\item[ii)] This result does not guarantee the
accuracy of the solution for the $N-n$ non conserved quantities.
In fact, since the very last step in the scheme is not a
projection towards the local equilibrium, a final layer effect
occurs. The use of stiffly accurate schemes (i.e.\ schemes for
which $a_{\nu j} = w_j, j=1,\ldots,\nu$) in the implicit step may
serve as a remedy to this problem.

\item[iii)] The theorem guarantees that in the stiff limit the
numerical scheme becomes the explicit RK scheme applied to the
equilibrium system, and therefore the order of accuracy of the
limiting scheme is greater or equal to the order of accuracy of
the original IMEX scheme.
\end{itemize}
\end{remarks}

When constructing numerical schemes for conservation laws, one has
to take a great care in order to avoid spurious numerical
oscillations arising near discontinuities of the solution. This is
avoided by a suitable choice of space discretization (see
Section~4) and time discretization.

Solution of scalar conservation equations, and equations with a
dissipative source have some norm that decreases in time. It would
be desirable that such property is maintained at a discrete level
by the numerical method. If $U^n$ represents a vector of solution
values (for example obtained from a method of lines approach in
solving (\ref{eq:fluid})) we recall the following \cite{SR}

\begin{definition}
A sequence $\{U^n\}$ is said to be {\em strongly stable} in a
given norm $||\cdot||$ provided that $||U^{n+1}||\leq ||U^n||$ for
all $n\geq 0$.
\end{definition}
The most commonly used norms are the $TV$-norm and the infinity norm. A
numerical scheme that maintains strong stability at discrete level is
called Strong Stability Preserving (SSP).

Here we review some basic facts about RK-SSP schemes. First, it
has been shown \cite{GST} under fairly general conditions that
high order SSP schemes are necessarily explicit. Second, observe
that a generic explicit RK scheme can be written as
\begin{eqnarray}
U^{(0)}  & = & U^n\nonumber \\
U^{(i)} & = & \sum_{k=0}^{i-1} (\alpha_{ik}U^{(k)} + \h \beta_{ik} L(U^{(k)}), \quad i=1,\ldots,\nu,\label{eq:RKSSP}\\
U^{n+1} & = & U^{(\nu)}\nonumber
\end{eqnarray}
where $\alpha_{ik}\geq0$ and $\alpha_{ik}=0$ only if
$\beta_{ik}=0$. This representation of RK schemes (which is not
unique) can be converted to a standard Butcher form in a
straightforward manner. Observe that for consistency, one has
$\sum_{k=0}^{i-1}\alpha_{ik} = 1$. It follows that if the scheme
can be written in the form \rf{eq:RKSSP} with non negative
coefficients $\beta_{ik}$ then it is a convex combination of
Forward Euler steps, with step sizes $\beta_{ik}/\alpha_{ik}\h$. A
consequence of this is that if Forward Euler is SSP for
$\h\leq\h^*$, then the RK scheme is also SSP for $\h\leq c\h^*$,
with $c=\min_{ik}(\alpha_{ik}/\beta_{ik})$ \cite{Shu-TVD,GST}.

The constant $c$ is a measure of the efficiency of the SSP-RK scheme,
therefore for the applications it is important to have $c$ as large as
possible. For a detailed description of optimal SSP schemes and their
properties see \cite{SR}.

By point iii) of the above remarks it follows that if
the explicit part of the IMEX scheme is SSP then, in the stiff
limit, we will obtain an SSP method for the limiting conservation
law. This property is essential to avoid spurious oscillations in
the limit scheme for equation \rf{eq:fluid}.

In this section we give the Butcher tableau for the new second and third order IMEX
schemes that satisfy the conditions of Theorem \ref{th:main}. In all these
schemes the implicit {\em tableau\/} corresponds to an L-stable scheme, that is
$w^T A^{-1}e=1$, $e$ being a vector whose components are all equal to 1
\cite{HW}, whereas the explicit tableau is SSP$k$, where $k$ denotes the order
of the SSP scheme. Several examples of asymptotically SSP schemes are reported
in tables (\ref{tb:PR222})-(\ref{tb:PR433}). We shall use the notation
SSP$k$$(s,\sigma, p)$, where the triplet $(s,\sigma, p)$ characterizes the
number $s$ of stages of the implicit scheme, the number $\sigma$ of stages of
the explicit scheme and the order $p$ of the IMEX scheme.

\begin{table}[h]
{\small
\[
    \begin{array}{c|cc}
      0 & 0  & 0       \\
      1 & 1  & 0       \\
      \hline
        & 1/2  & 1/2
    \end{array}  \qquad
    \begin{array}{c|cccc}
      \gamma & \gamma        & 0              \\
      1-\gamma   & 1-2\gamma & \gamma \\
      \hline
          & 1/2   &  1/2
    \end{array} \qquad
      \gamma = 1-\frac1{\sqrt{2}}
\]
\caption{Tableau for the explicit (left) implicit (right)
IMEX-SSP2(2,2,2) L-stable scheme} \label{tb:PR222} }
\end{table}

\begin{table}[h]
{\small
\[
    \begin{array}{c|ccc}
      0 & 0  & 0  & 0     \\
      0 & 0  & 0  & 0     \\
      1 & 0  & 1  & 0     \\      \hline
        & 0  & 1/2 & 1/2
    \end{array}  \qquad
    \begin{array}{c|cccc}
      1/2 & 1/2   & 0 & 0              \\
      0 & -1/2 & 1/2  & 0              \\
         1 & 0 & 1/2 &  1/2 \\
      \hline
          & 0 & 1/2 &  1/2
    \end{array}
\]
\caption{Tableau for the explicit (left) implicit (right)
IMEX-SSP2(3,2,2) stiffly accurate scheme} \label{tb:PR322}
}
\end{table}

\begin{table}[h]
{\small
\[
    \begin{array}{c|ccc}
      0 & 0  & 0  & 0     \\
      1/2 & 1/2  & 0  & 0     \\
      1 & 1/2  & 1/2  & 0     \\      \hline
        & 1/3  & 1/3 & 1/3
    \end{array}  \qquad
    \begin{array}{c|cccc}
      1/4 & 1/4   & 0 & 0              \\
      1/4 & 0 & 1/4  & 0              \\
         1 & 1/3 & 1/3 &  1/3 \\
      \hline
          & 1/3 & 1/3 &  1/3
    \end{array}
\]
\caption{Tableau for the explicit (left) implicit (right)
IMEX-SSP2(3,3,2) stiffly accurate scheme} \label{tb:PR332}
}
\end{table}

\begin{table}[h]
{\small
\[
    \begin{array}{c|ccc}
      0 & 0  & 0  & 0     \\
      1 & 1  & 0  & 0     \\
      1/2 & 1/4  & 1/4  & 0     \\      \hline
        & 1/6  & 1/6 & 2/3
    \end{array}  \qquad
    \begin{array}{c|ccc}
      \gamma & \gamma & 0 & 0              \\
      1-\gamma & 1-2\gamma & \gamma  & 0     \\
      1/2 & 1/2-\gamma & 0 &  \gamma \\
      \hline
          & 1/6 & 1/6 &  2/3
    \end{array} \qquad
      \gamma = 1-\frac1{\sqrt{2}}
\]
\caption{Tableau for the explicit (left) implicit (right)
IMEX-SSP3(3,3,2) L-stable scheme} \label{tb:PR332b} }
\end{table}

\begin{table}[h]
{\small
\[
    \begin{array}{c|cccc}
      0 & 0 & 0  & 0  & 0     \\
      0 & 0 & 0  & 0  & 0     \\
      1 & 0 & 1  & 0  & 0     \\
      1/2 & 0 & 1/4  & 1/4  & 0     \\      \hline
        & 0 & 1/6  & 1/6 & 2/3
    \end{array}  \qquad
    \begin{array}{c|cccc}
      \alpha & \alpha & 0 & 0 & 0              \\
           0 & -\alpha & \alpha  & 0  & 0            \\
           1 & 0 & 1-\alpha &  \alpha & 0 \\
          1/2& \beta & \eta & 1/2-\beta-\eta-\alpha & \alpha\\
      \hline
          & 0 & 1/6 & 1/6 &  2/3
    \end{array}
\]
\[
\alpha = 0.24169426078821,\quad \beta=0.06042356519705\quad
\eta=0.12915286960590
\]
\caption{Tableau for the explicit (left) implicit (right)
IMEX-SSP3(4,3,3) L-stable scheme} \label{tb:PR433} }
\end{table}

\section{IMEX-WENO schemes}
For simplicity we consider the case of the single scalar equation
\begin{equation}
  u_t + f(u)_x = \frac1{\e} \r(u).
  \label{eq:scal}
\end{equation}
We have to distinguish between schemes based on cell averages (finite volume
approach, widely used for conservation laws) and schemes based on point values
(finite difference approach).

Let {$\hh$} and {$\Delta t$} be the mesh widths. We introduce the
grid points $$ x_j=j\hh, \quad x_{j+1/2}=x_j+\frac12{\hh}, \quad
j=\ldots,-2,-1,0,1,2,\ldots$$ and use the standard notations
$$ u_j^n=u(x_j,t^n),\qquad
\bar{u}_{j}^{n}=\frac1{\hh}\int_{x_{j-1/2}}^{x_{j+1/2}}
u(x,t^{n})\,dx.
$$

\subsection{Finite volumes}

Integrating equation (\ref{eq:scal}) on {$I_j=[x_{j-1/2},x_{j+1/2}]$} and
dividing by {$\Delta x$} we obtain
\begin{equation}
   \frac{d\bar{u}_j}{dt} = -\frac1{\hh} [f(u(x_{j+1/2},t))-f(u(x_{j-1/2},t))] +\frac1{\e} \overline{r(u)}_j.
   \label{eq:FV}
\end{equation}
In order to convert this expression into a numerical scheme, one has to
approximate the right hand side with a function of the cell averages
$\{\bar{u}(t)\}_j$, which are the basic unknowns of the problem.

The first step is to perform a reconstruction of the unknown function
$u(x,t)$ by a piecewise polynomial, starting from cell averages
$\bar{u}_j(t)$. Given $\{\bar{u}_{j}\}$, compute a piecewise polynomial
reconstruction
$$
\Rz(x)=\sum_{j} P_{j}(x)\chi_{j}(x),
$$
where $P_j(x)$ is a polynomial satisfying some accuracy and non
oscillatory property, and $\chi_{j}(x)$ is the indicator function of
cell $I_j$. For first order schemes, the  reconstruction is piecewise
constant, while second order schemes can be obtained by a piecewise
linear reconstruction. Higher order schemes are obtained by higher
order polynomials.

The reconstruction step is very delicate, because the function $u(x,t)$
may be discontinuous. To this goal, suitable techniques, such as
Weighted Essentially Non Oscillatory (WENO), have been developed. The
reader can consult the review by Shu \cite{Shu-CIME} and references
therein for a more detailed description of high order non oscillatory
reconstructions.

The flux function at the edge of the cell can be computed by using a suitable
numerical flux function, consistent with the analytical flux,
\[
  f(u(x_{j+1/2},t)) \approx F(u_{j+1/2}^{-},u_{j+1/2}^{+}),
\]
where the quantities $u_{j+1/2}^\pm$ are obtained from the
reconstruction, as $u_{j+1/2}^\pm = \lim_{x\to x_{j+1/2}^\pm}\Rz(x)$.

Example of numerical flux functions are the Godunov flux
\[
  F(u_{j+1/2}^-,u_{j+1/2}^+) = f(u^*(u_{j+1/2}^-,u_{j+1/2}^+)),
\]
where $u^*(u_{j+1/2}^-,u_{j+1/2}^+)$ is the solution of the Riemann problem at
$x_{j+1/2}$, corresponding to the states $u_{j+1/2}^-$ and $u_{j+1/2}^+$, and
the Local Lax Friedrichs flux (also known as Rusanov flux),
\[
  F(u_{j+1/2}^-,u_{j+1/2}^+) =
  \frac12[(f(u_{j+1/2}^-)+f(u_{j+1/2}^+))-\alpha(u_{j+1/2}^+-u_{j+1/2}^-)],
\]
where $\alpha=\max_{w}|f^{\prime}(w)|$, and the maximum is taken over the
relevant range of $w$.

The two examples constitute two extreme cases of numerical fluxes:
the Godunov flux is the most accurate and the one that produces
the best results for a given grid size, but it is very expensive,
since it requires the solution to the Riemann problem. Local
Lax-Friedrichs flux, on the other hand, is less accurate, but much
cheaper. This latter is the numerical flux that has been used
throughout all the calculations performed for this paper. The
difference in resolution provided by the various numerical fluxes
becomes less relevant with the increase in the order of accuracy
of the method.

The right hand side of Eq.(\ref{eq:FV}) contains the average of the source term
{$\overline{r(u)}$} instead of the source term evaluated at the average of
{$u$}, {$r(\bar{u})$}. The two quantities agree within second order accuracy
\[
\overline{r(u)}_j=r(\bar{u}_j)+O(\hh^2).
\]
This approximation can be used to construct schemes up to second order.

First order (in space) semidiscrete schemes can be obtained using the numerical
flux function $F(\bar{u}_{j},\bar{u}_{j+1})$ in place of $f(u(x_{j+1/2},t))$,
\begin{equation}
\frac{d\bar{u}_j}{dt}=-\frac{F(\bar{u}_{j},\bar{u}_{j+1})-
F(\bar{u}_{j-1},\bar{u}_{j})}{\hh}+\frac1{\e} {r(\bar{u}_j)}.
\label{eq:secondFV}
\end{equation}

Second order schemes are obtained by using a piecewise linear reconstruction in
each cell, and evaluating the numerical flux on the two sides of the interface
$$
\frac{d\bar{u}}{dt} =-\frac{F(u_{j+1/2}^{-},u_{j+1/2}^{+})-
F(u_{j-1/2}^{-},u_{j-1/2}^{+})}{\hh} + \frac{1}{\e}r(\bar {u}_j)
$$
The quantities at cell edges are computed by piecewise linear
reconstruction. For example,
\[
   u_{j+1/2}^- = \bar{u}_j + \frac{\hh}{2}u'_j
\]
where the slope $u'_j$ is a first order approximation of the space derivative
of $u(x,t)$, and can be computed by suitable slope limiters (see, for example,
\cite{leveque:numerical-methods} for a discussion on TVD slope limiters.)

For schemes of order higher than second, a suitable quadrature
formula is required to approximate $\overline{g(u)}_j$. For
example, for third and fourth order schemes, one can use Simpson's
rule
\[
  \overline{r(u)}_j\approx \frac16(r(u_{j-1/2}^+) + 4 r(u_j) +
  r(u_{j+1/2}^-))
\]
where the pointwise values $u_{j-1/2}^+$, $u_j$, $u_{j+1/2}^-$ are
obtained from the reconstruction.

For a general problem, this has the effect that the source term couples
the cell averages of different cells, thus making almost impractical
the use of finite volume methods for high order schemes applied to
stiff sources.

Note, however, that in many relevant cases of hyperbolic systems with
relaxation the implicit step, thanks to the conservation properties of
the system, can be explicitly solved, and finite volume methods can be
successfully used. We mention here all relaxation approximation of
Jin-Xin type \cite{JX}, some simple discrete velocity models, such as
Carlemann and Broadwell models, monoatomic gas in Extended
Thermodynamics, semiconductor models, and shallow water equations.

\subsection{Finite differences}
In a finite difference scheme the basic unknown is the pointwise value of the
function, rather than its cell average. Osher and Shu observed that it is
possible to write a finite difference scheme in conservative form
\cite{Shu-OsherII}. Let us consider the equation
$$
\devp{u}{t}+\devp{f}{x}=\frac{1}{\e}r(u),
$$
and write
$$
\devp{f}{x}(u(x,t))=\frac{\hat{f}(u(x+\frac{\hh}{2},t))-\hat{f}(u(x-\frac{\hh}{2},t))}{\hh}.
$$
The relation between $f$ and $\hat{f}$ is the
following. Let us consider the sliding average operator
$$
   \bar{u}(x,t) = <u>_x \equiv \frac{1}{\hh}
   \int_{x-\frac{\hh}{2}}^{x+\frac{\hh}{2}}{u(\xi,t)\,d\xi}.
$$
Differentiating with respect to $x$ one has
$$
\devp{\bar{u}}{x}=\frac{1}{\hh}(u(x+\frac{\hh}{2},t)-u(x-\frac{\hh}{2},t)).
$$
Therefore the relation between $f$ and $\hat{f}$ is the same that
exists between $\bar{u}(x,t)$ and $u(x,t)$, namely, flux function $f$
is the cell average of the function $\hat{f}$. This also suggests
a way to compute the flux function. The technique that is used to
compute pointwise values of $u(x,t)$ at the edge of the cell from
cell averages of $u$ can be used to compute
$\hat{f}(u(x_{j+1/2},t))$ from $f(u(x_{j},t))$. This means that in the
finite difference method it is the flux function which is computed
at $x_{j}$ and then reconstructed at $x_{j+1/2}$. But the
reconstruction at $x_{j+1/2}$ may be discontinuous. Which value
should one use? A general answer to this question can be given if
one considers flux functions that can be split
\begin{equation}
\label{split} f(u)=f^{+}(u)+f^{-}(u)\sv
\end{equation}
with the condition that
\begin{equation}
   \der{f^{+}(u)}{u}\geq 0\sv \quad \der{f^{-}(u)}{u}\leq 0.
   \label{eq:monotone}
\end{equation}
There is a close analogy between flux splitting and numerical flux functions.
In fact, if a flux can be split as \eq{split}, then
$$
F(a,b)=f^{+}(a)+f^{-}(b)
$$
will define a monotone consistent flux, provided condition
(\ref{eq:monotone}) is satisfied. Together with non oscillatory
reconstructions (such as WENO) and SSP time discretization, the
monotonicity condition will ensure that the overall scheme will not
produce spurious numerical oscillations (see, for example,
\cite{leveque:numerical-methods, Shu-CIME}.)

This is the case, for example, of the local Lax-Friedrichs flux. A
finite difference scheme therefore takes the following form
$$
\der{u_{j}}{t}=-\frac{1}{\hh}[\hat{F}_{j+1/2}-\hat{F}_{j-1/2}] + \frac1\e g(u_j),
$$
$$
\hat{F}_{j+1/2}=\hat{f}^{+}(u_{j+1/2}^{-})+\hat{f}^{-}(u_{j+1/2}^{+});
$$
$\hat{f}^{+}(u_{j+1/2}^{-})$ is obtained by
\begin{itemize}
\item computing $f^{+}(u_{l})$ and interpret it as cell average of
$\hat{f}^{+}$;
\item performing pointwise reconstruction of $\hat{f}^{+}$ in cell $j$, and
evaluate it in $x_{j+1/2}$.
\end{itemize}
$\hat{f}^{-}(u_{j+1/2}^{+})$ is obtained by
\begin{itemize}
\item computing $f^{-}(u_{l})$, interpret as cell average of
$\hat{f}^{-}$;
\item performing pointwise reconstruction of $\hat{f}^{-}$ in cell $j+1$, and
evaluate it in $x_{j+1/2}$.
\end{itemize}
A detailed account on high order finite difference schemes can be found in
\cite{Shu-CIME}.

At variance with finite volume schemes, since the source is evaluated
pointwise, finite difference schemes do not couple the cells. This is a
big advantage in the cases in which the implicit step can not be
explicitly solved.

\begin{remarks}
\begin{itemize}
\item[i)] Finite difference can be used only with uniform (or smoothly
varying) mesh. In this respect finite volume are more flexible,
since they can be used even with unstructured grids.

\item[ii)] The treatment presented for the scalar equation can be
extended to systems, with a minor change of notation. In
particular, the parameter $\alpha$ appearing in the local
Lax-Friedrichs flux will be computed using the spectral radius of
the Jacobian matrix.

\item[iii)] For schemes applied to systems, better results are usually
obtained if one uses characteristic variables rather than
conservative variables in the reconstruction step. The use of
conservative variables may result in the appearance of small
spurious oscillations in the numerical solution. For a treatment
of this effect see, for example, \cite{Shu-central}.
\end{itemize}

\end{remarks}

\section{Numerical tests}

In this section we investigate numerically the convergence rate
and the zero relaxation limit behavior of the schemes. To this aim
we apply the IMEX-WENO schemes to the Broadwell equations of
rarefied gas dynamics \cite{CJR, Jin-RK, LRR, Pa}. In all the
computations presented in this paper we used finite difference
WENO schemes with Lax-Friedrichs flux and conservative variables.
Of course the sharpness of the resolution of the numerical results
can be improved using a less dissipative flux.

As a comparison, together with the new IMEX-SSP schemes, we have
considered the second order ARS(2,2,2) method presented in \cite{ARS},
which was the one recommended by the authors. For this scheme, since
$c_1=0$, there will be a degradation of accuracy due to the initial
layer.

In all figures, the value of $N$ represents the number of grid points
in space.

The kinetic model is characterized by a hyperbolic system with
relaxation of the form (\ref{eq:hsr}) with
\[
U=(\rho,\,m,\,z), \qquad F(U)=(m,\,z,\,m), \qquad
R(U)=\left(0,\,0,\,\frac12{(\rho^2+m^2-2\rho z)}\right).
\]
Here $\e$ represents the mean free path of particles. The only
conserved quantities are the density $\rho$ and the momentum $m$.

In the fluid-dynamic limit $\e \to 0$ we have \be
z=z_E\equiv\frac{\rho^2+m^2}{2\rho}, \ee and the Broadwell system
is well approximated by the reduced system (\ref{eq:fluid}) with
\[
u=(\rho,\, \rho v),\qquad G(u)=\left(\rho v,\, \frac12(\rho + \rho v^2)\right),
\qquad v = \frac{m}{\rho},
\]
which represents the corresponding Euler equations of
fluid-dynamics.

\subsubsection*{Convergence rates}

\begin{figure}[htb]
\centering{
\includegraphics[scale=0.41]{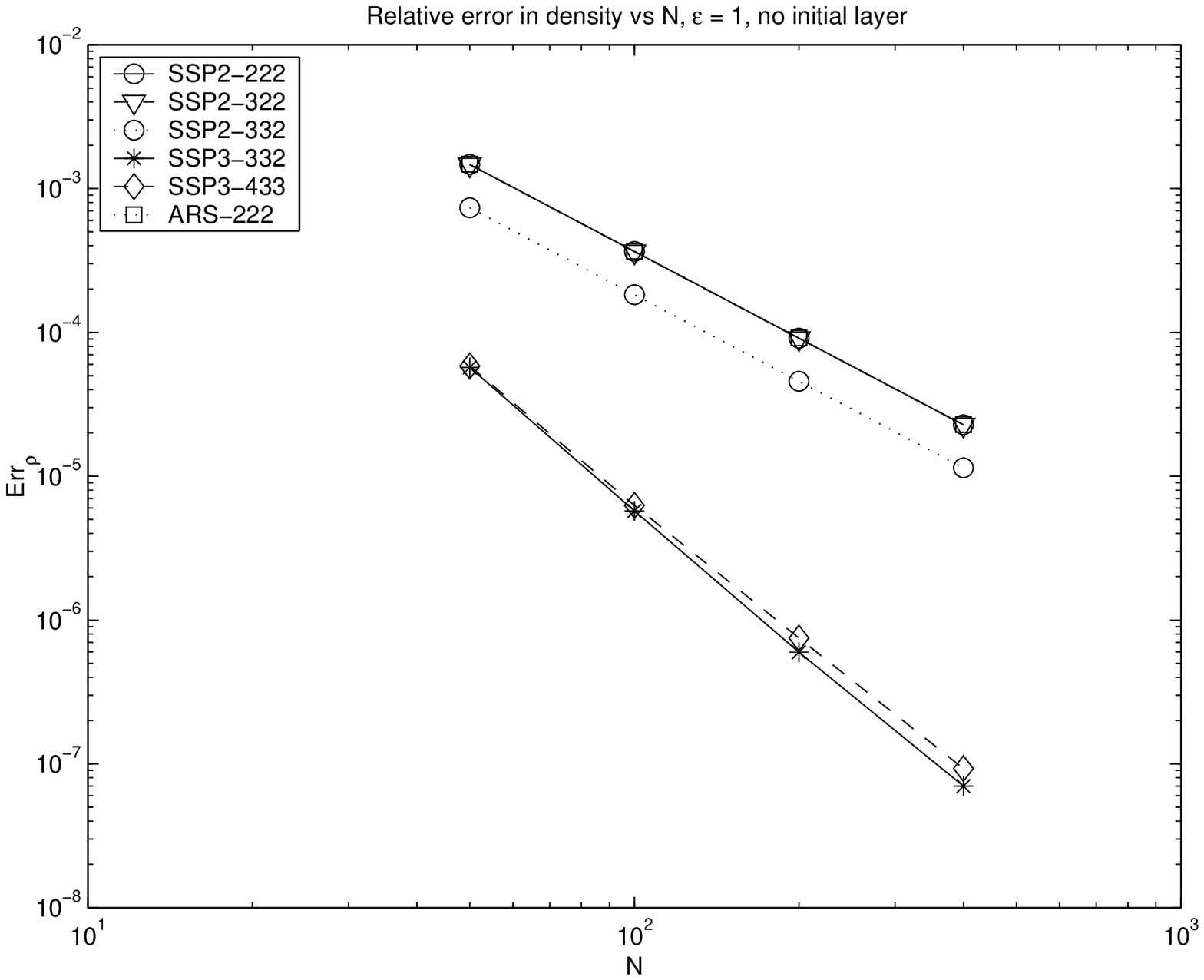}
\hskip .5cm
\includegraphics[scale=0.41]{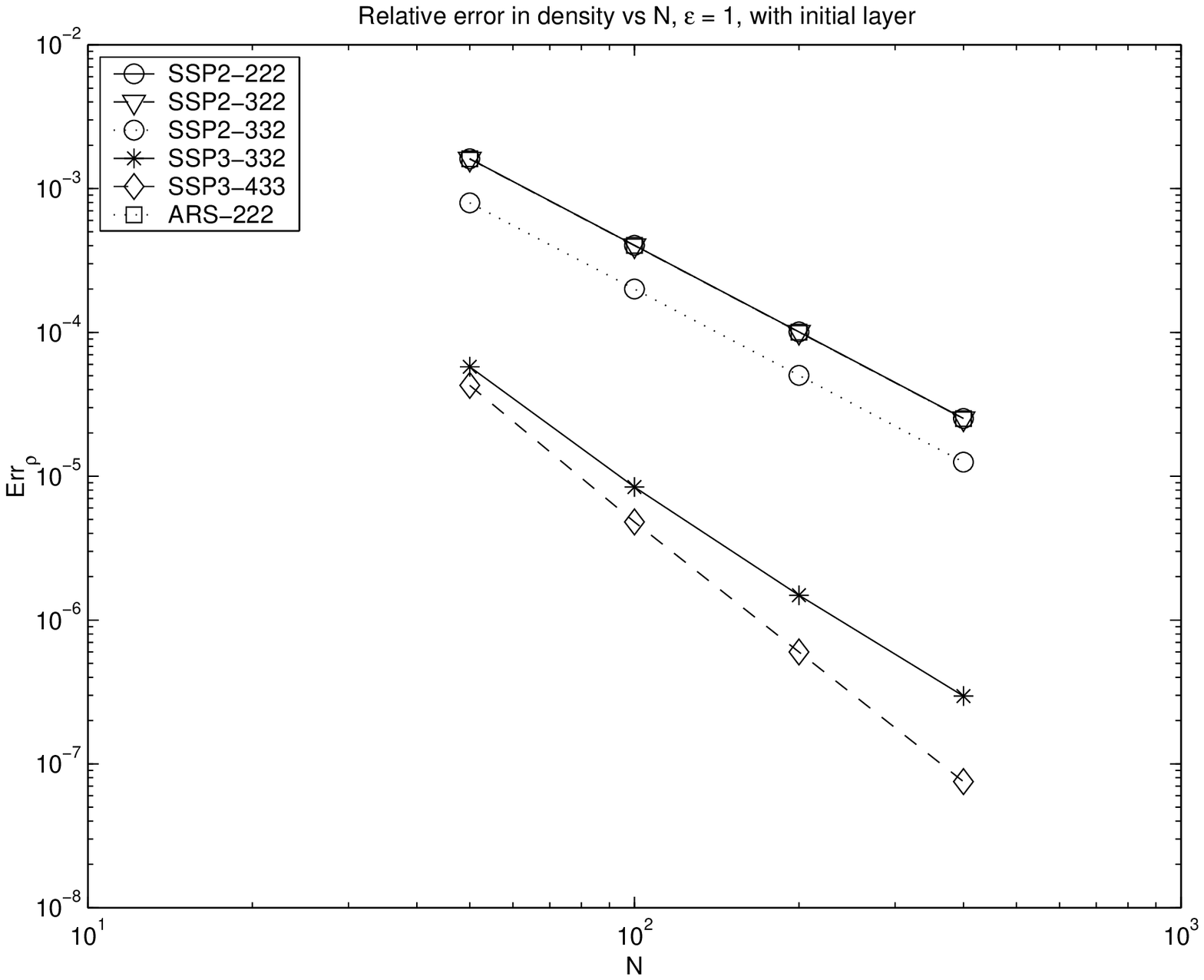}
\includegraphics[scale=0.41]{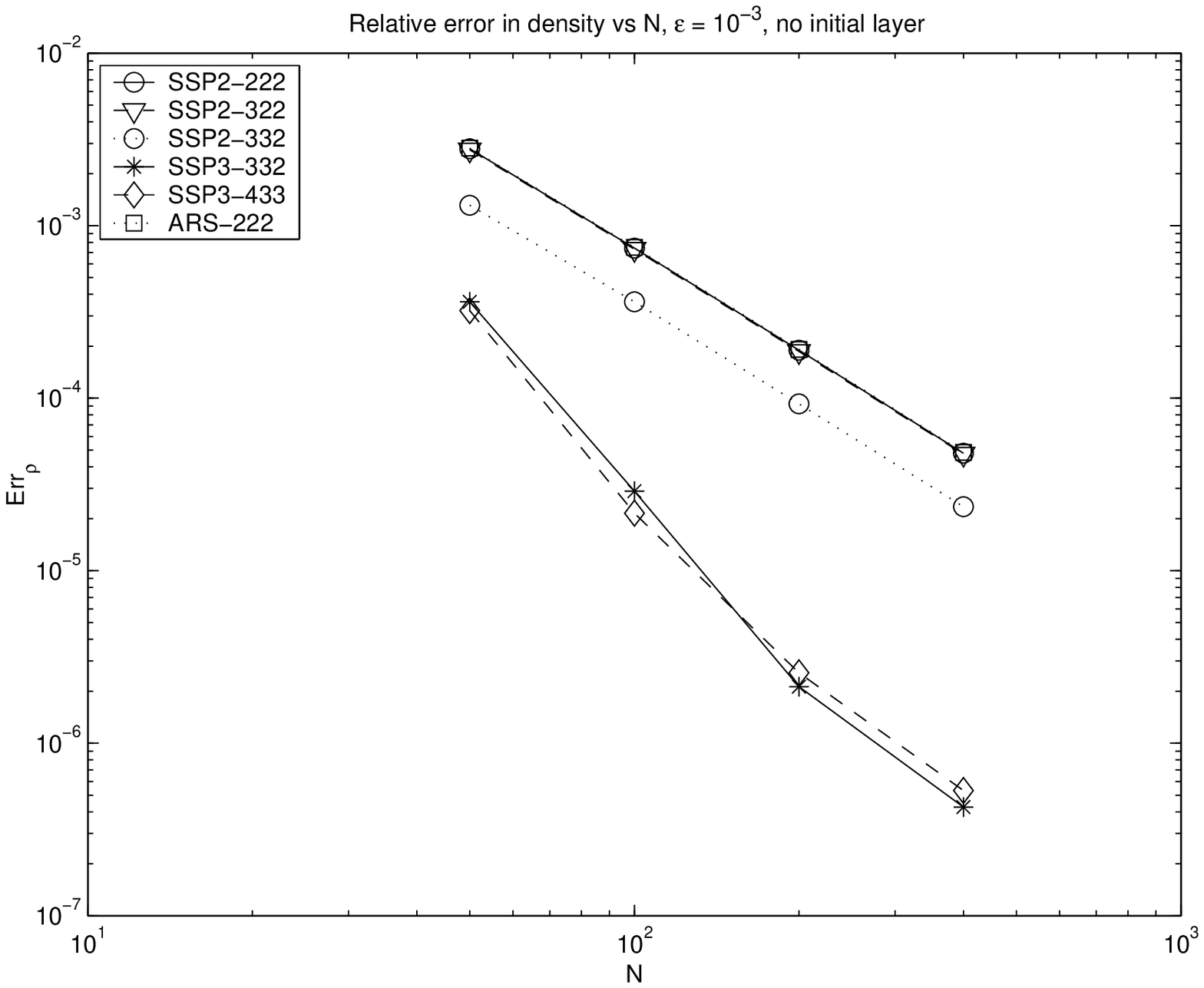}
\hskip .5cm
\includegraphics[scale=0.41]{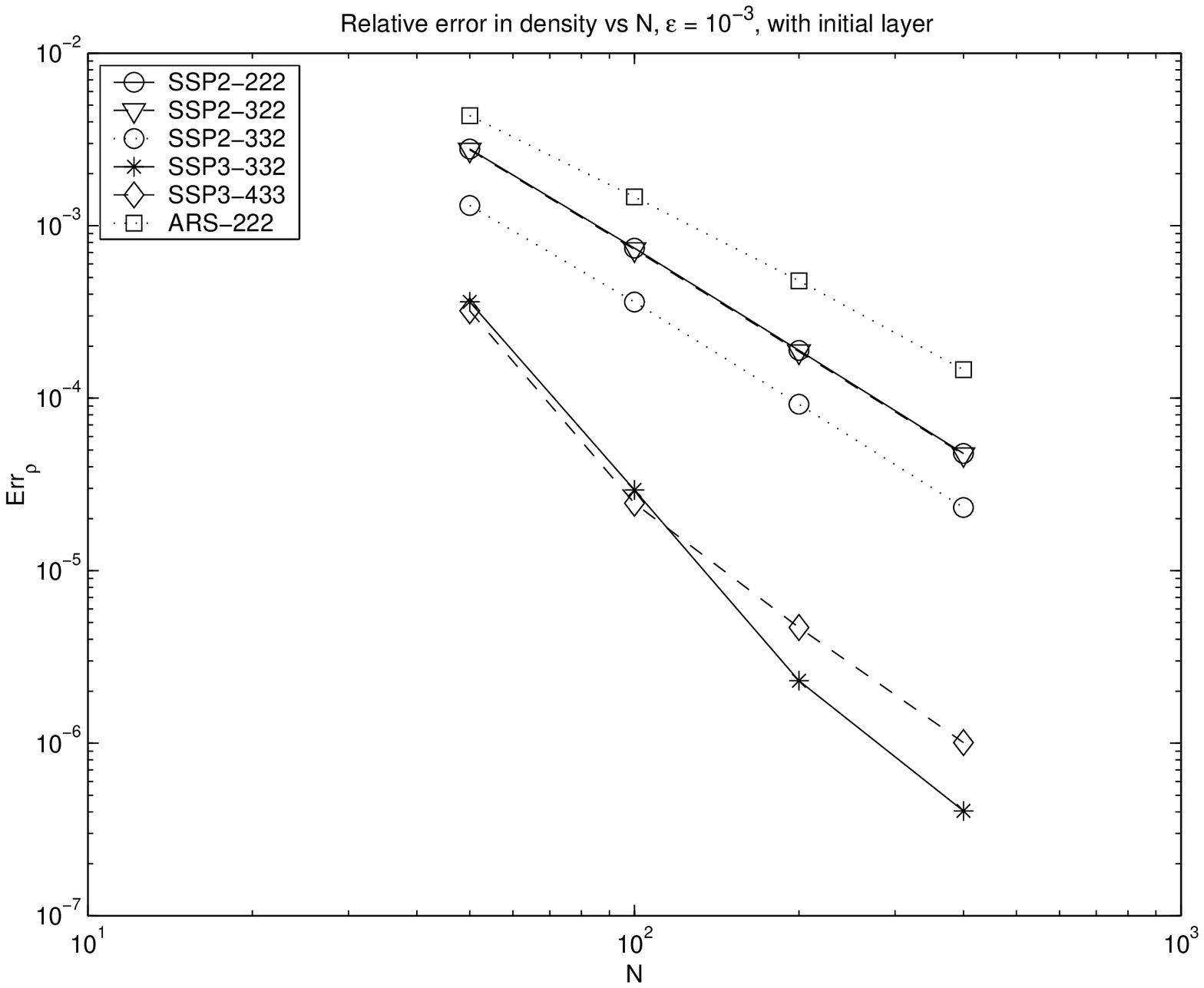}
\includegraphics[scale=0.41]{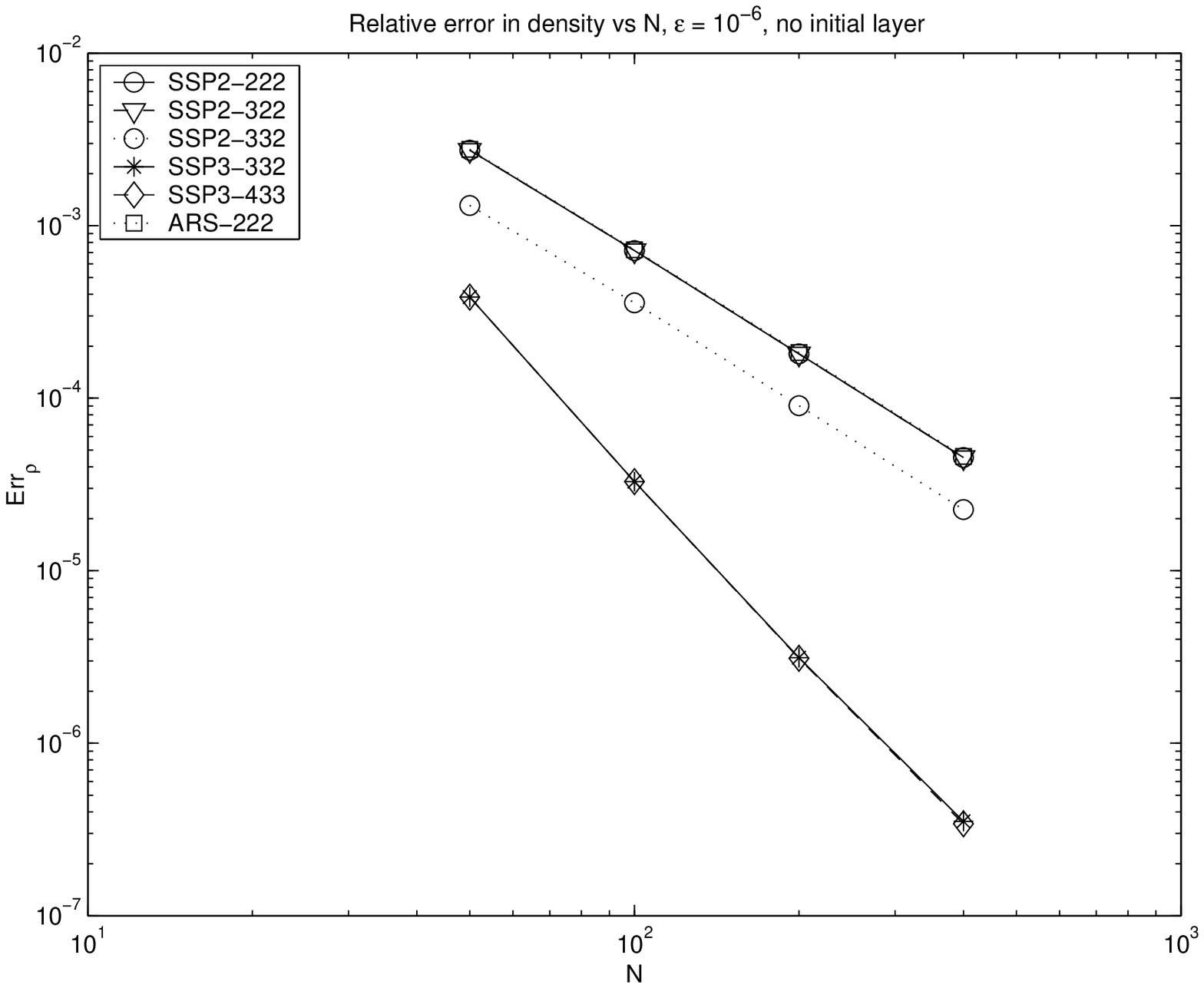}
\hskip .5cm
\includegraphics[scale=0.41]{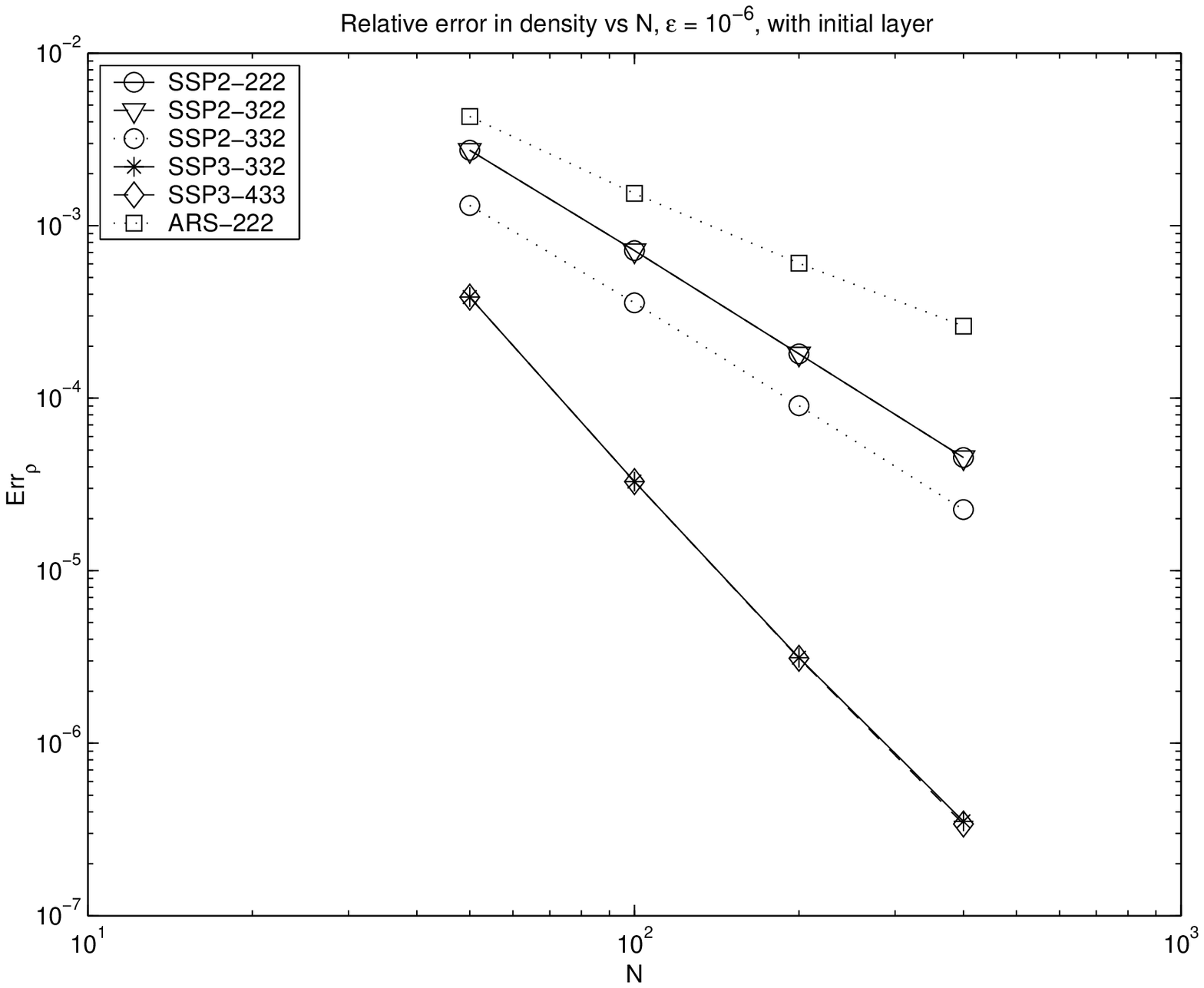}
} \vskip .5cm \caption{Relative errors for density $\rho$ in the
Broadwell equations with initial data (\ref{eq:conv}). Left column
$a_z=1.0$ (no initial layer), right column $a_z=0.2$ (initial
layer). From top to bottom, $\e=1.0,\, 10^{-3},\, 10^{-6}$.}
\label{fig:c}
\end{figure}

We have considered a periodic smooth solution with initial data as
in \cite{CJR, LRR} given by
\be
   \rho(x,0)=1+a_{\rho}\sin\frac{2\pi
   x}{L}, \quad v(x,0)=\frac{1}{2}+a_{v}\sin\frac{2\pi x}{L},\quad
   z(x,0)=a_{z}\frac{\rho(x,0)^2+m(x,0)^2}{2\rho(x,0)}.
\ee
In our computations we used the parameters
\begin{equation}
   a_\rho=0.3,\quad a_v = 0.1,
   \quad a_z= 1.0 \mbox{ (no initial layer) and }  a_z = 0.2 \mbox{ (initial layer)}, \quad L=20,
   \label{eq:conv}
\end{equation}
and we integrate the equations for $t \in [0,30]$. A Courant number $\Delta
t/\Delta x = 0.6$ has been used. The plots of the relative error are given in
Figure \ref{fig:c}.

Notice how, in absence of initial an layer, all schemes tested have the
prescribed order of accuracy both in the non stiff and in the stiff
limit, with some degradation of the accuracy at intermediate regimes.
As expected, scheme ARS(2,2,2), shows a degradation of the accuracy
when an initial layer is present.

\begin{figure}[htb]
\centering{
\includegraphics[scale=0.41]{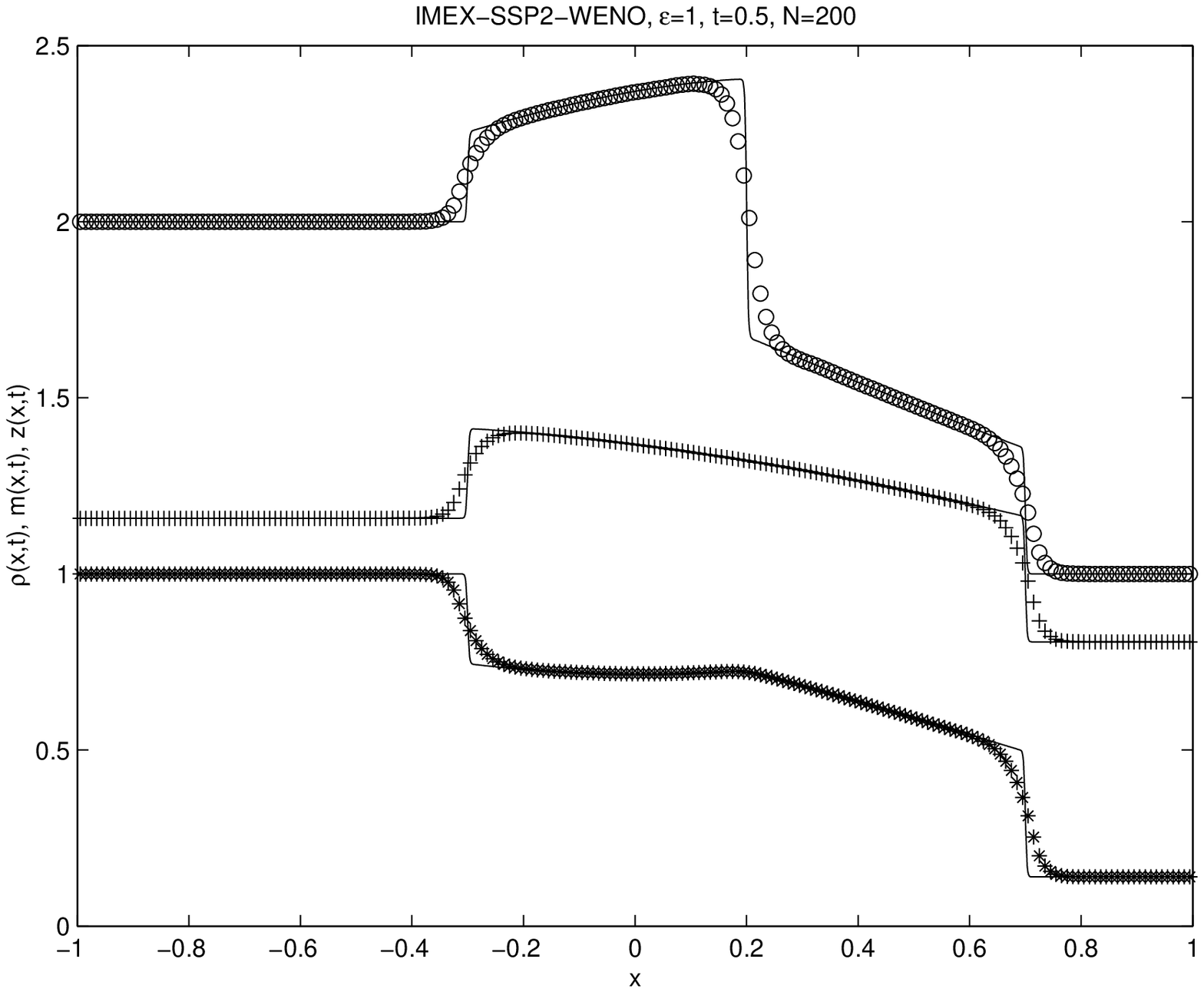}
\hskip .5cm
\includegraphics[scale=0.41]{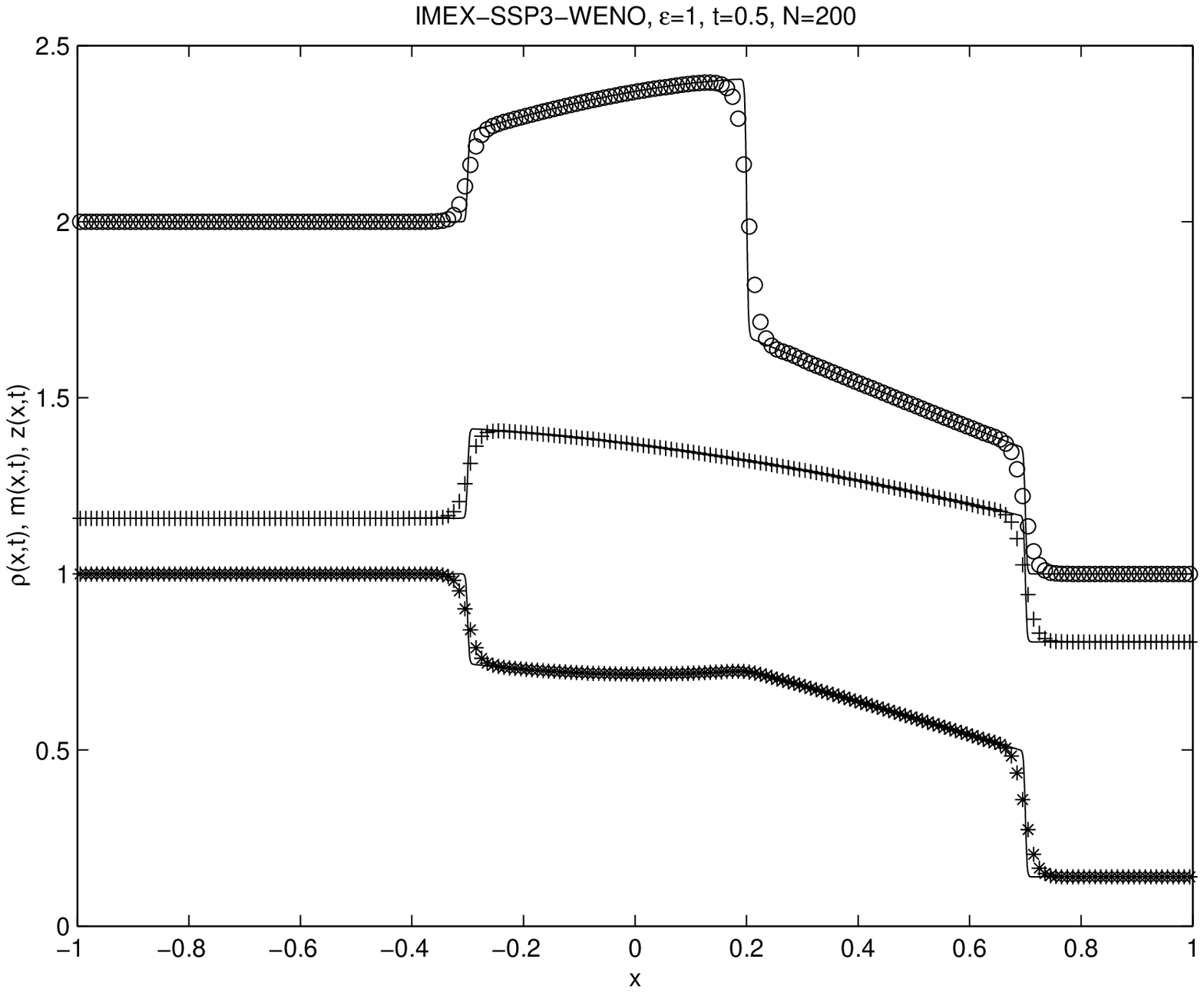}
\includegraphics[scale=0.41]{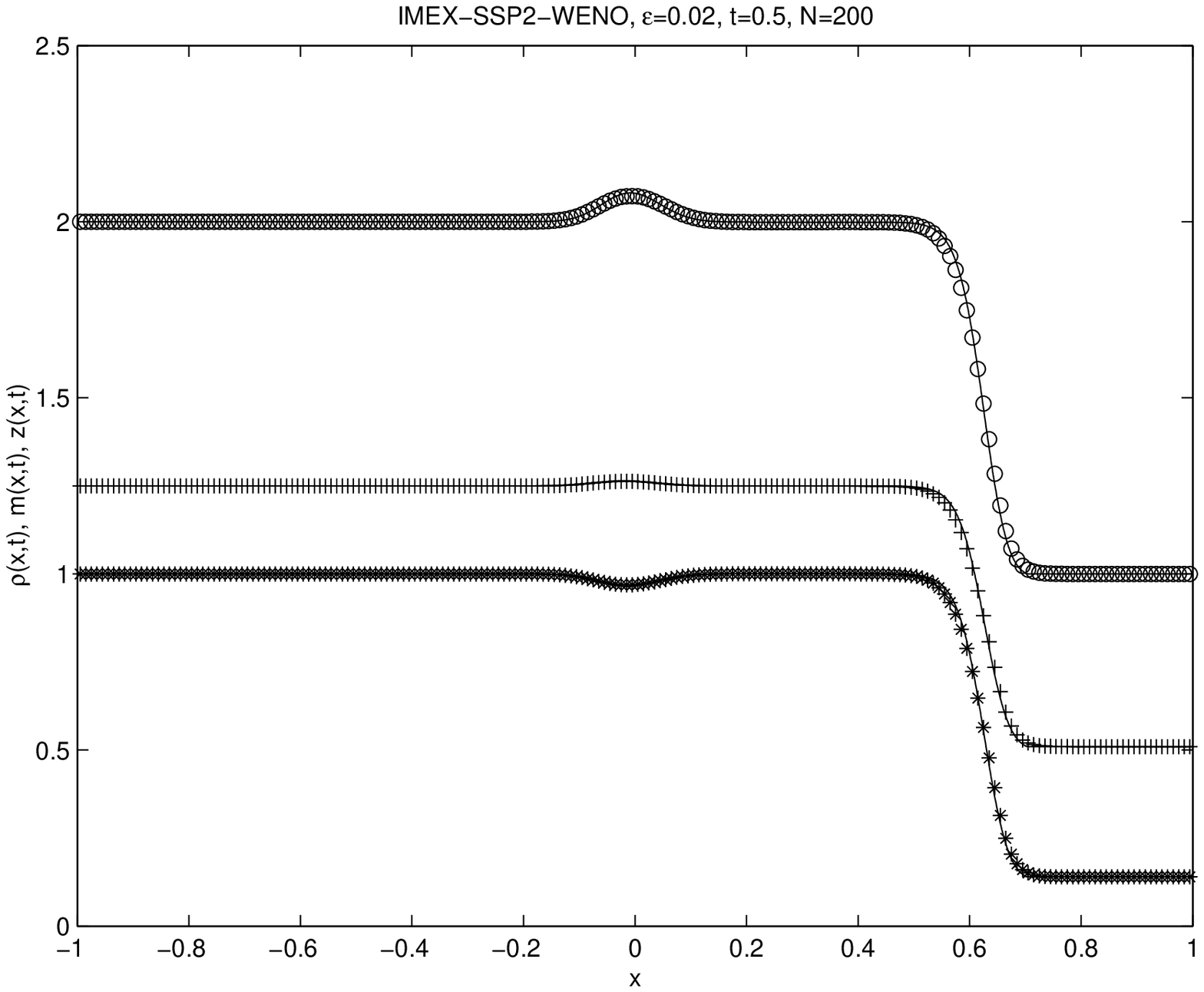}
\hskip .5cm
\includegraphics[scale=0.41]{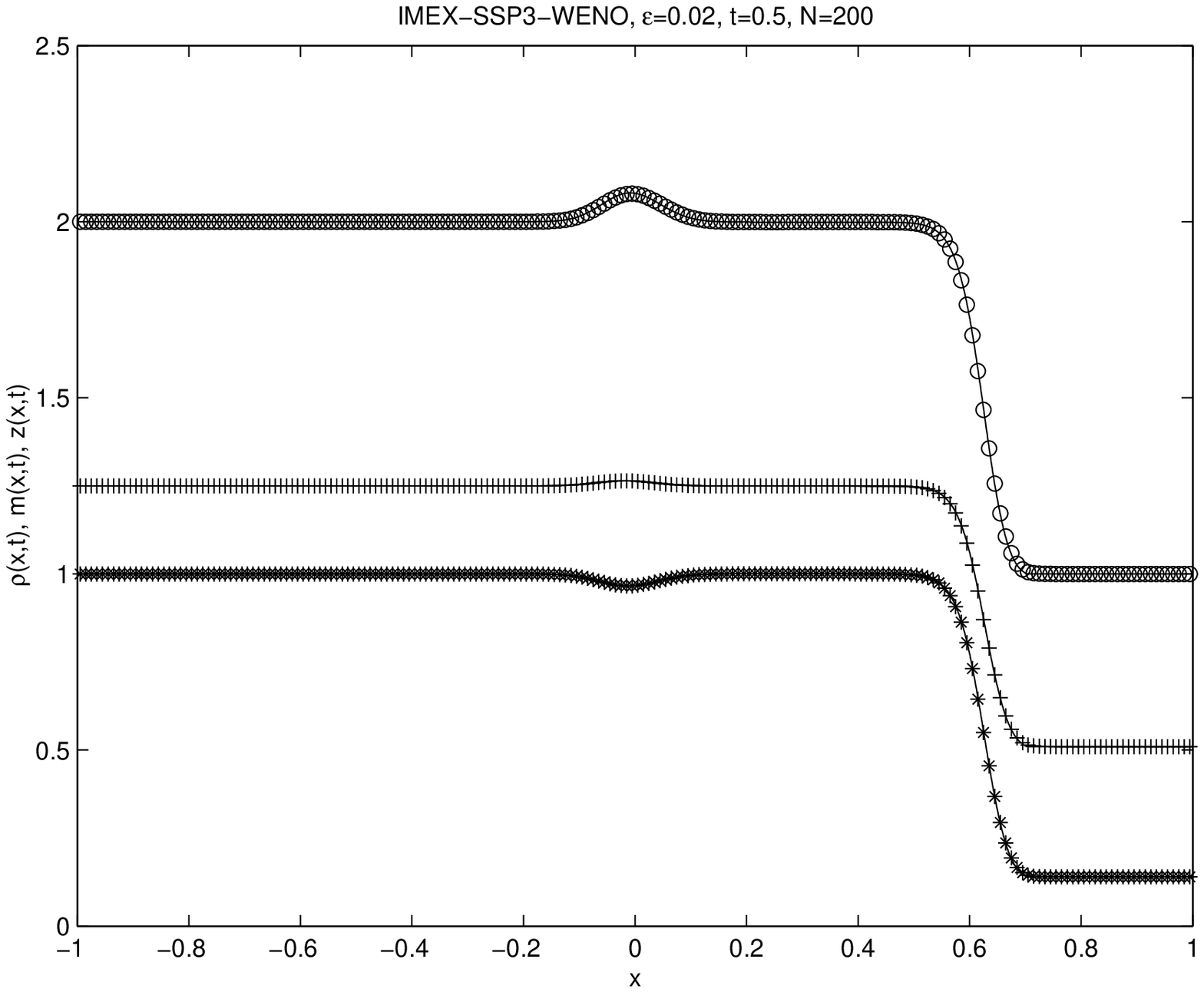}
\includegraphics[scale=0.41]{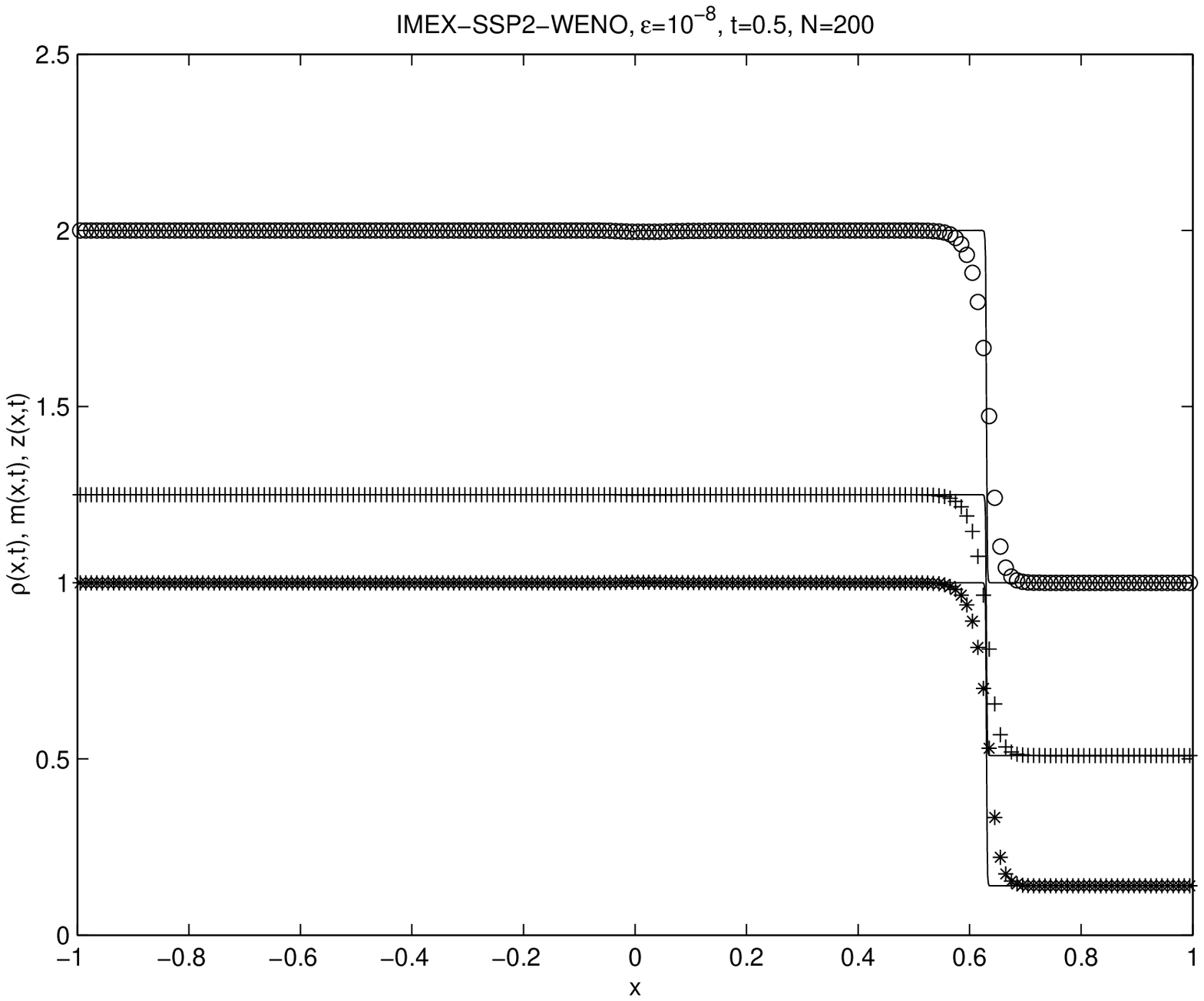}
\hskip .5cm
\includegraphics[scale=0.41]{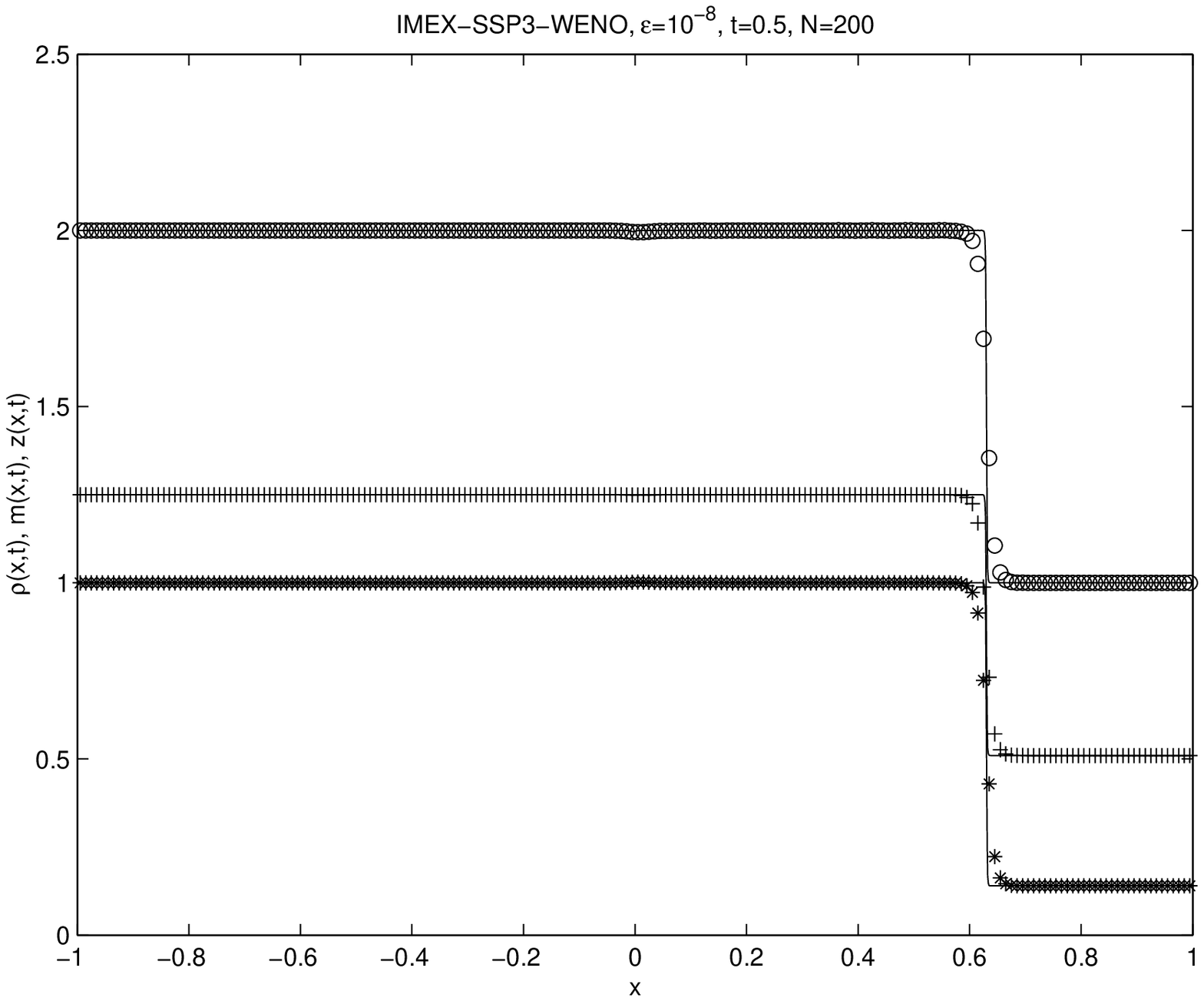}
} \vskip .5cm \caption{Numerical solution of the Broadwell
equations with initial data (\ref{eq:rim1}) for $\rho$($\circ$),
$m$($\ast$) and $z$($+$) at time $t=0.5$. Left column
IMEX-SSP2-WENO scheme, right column IMEX-SSP3-WENO scheme. From
top to bottom, $\e=1.0,\, 0.02,\, 10^{-8}$.} \label{fig:1}
\end{figure}

\begin{figure}[htb]
\centering{
\includegraphics[scale=0.41]{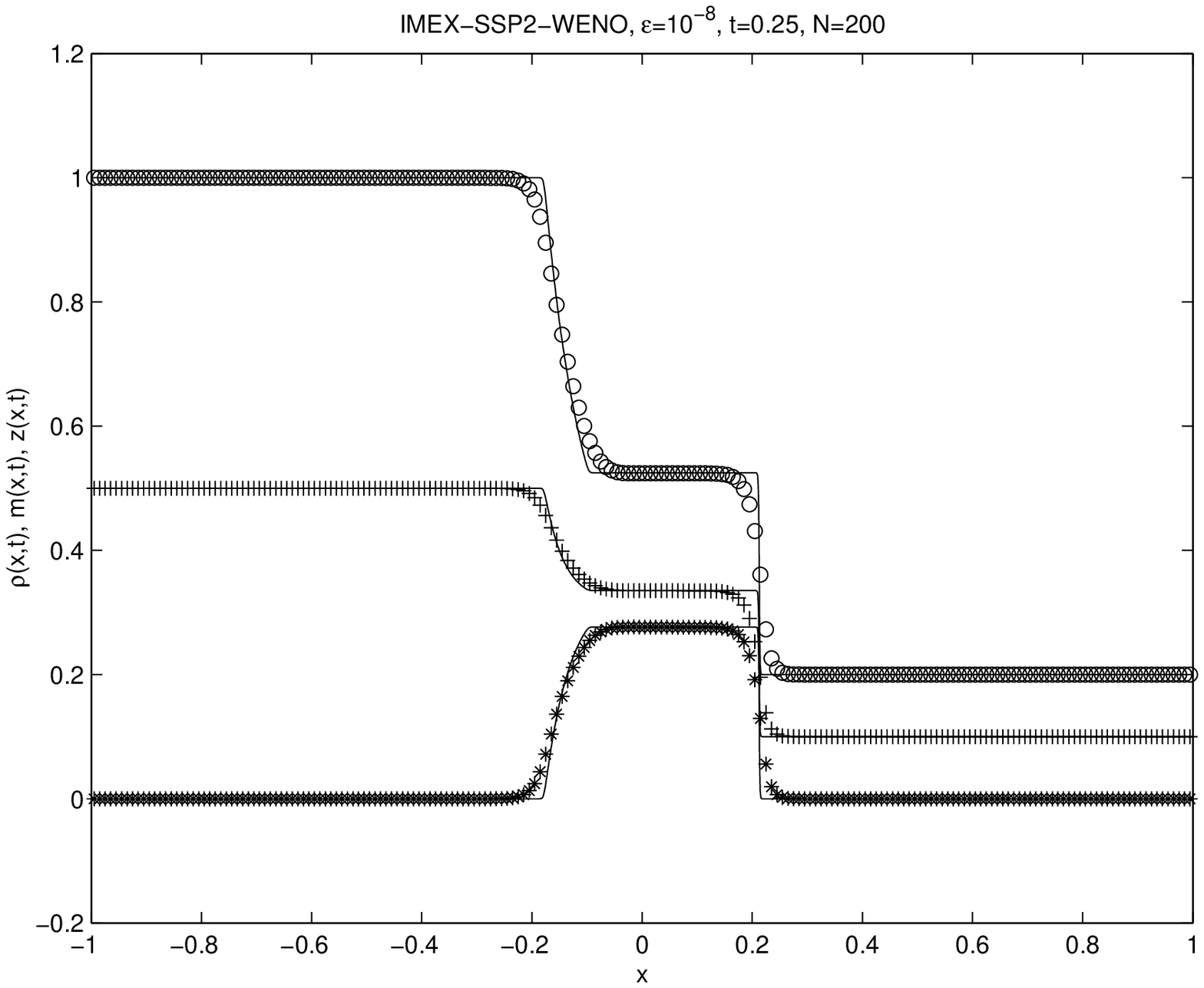}
\hskip .5cm
\includegraphics[scale=0.41]{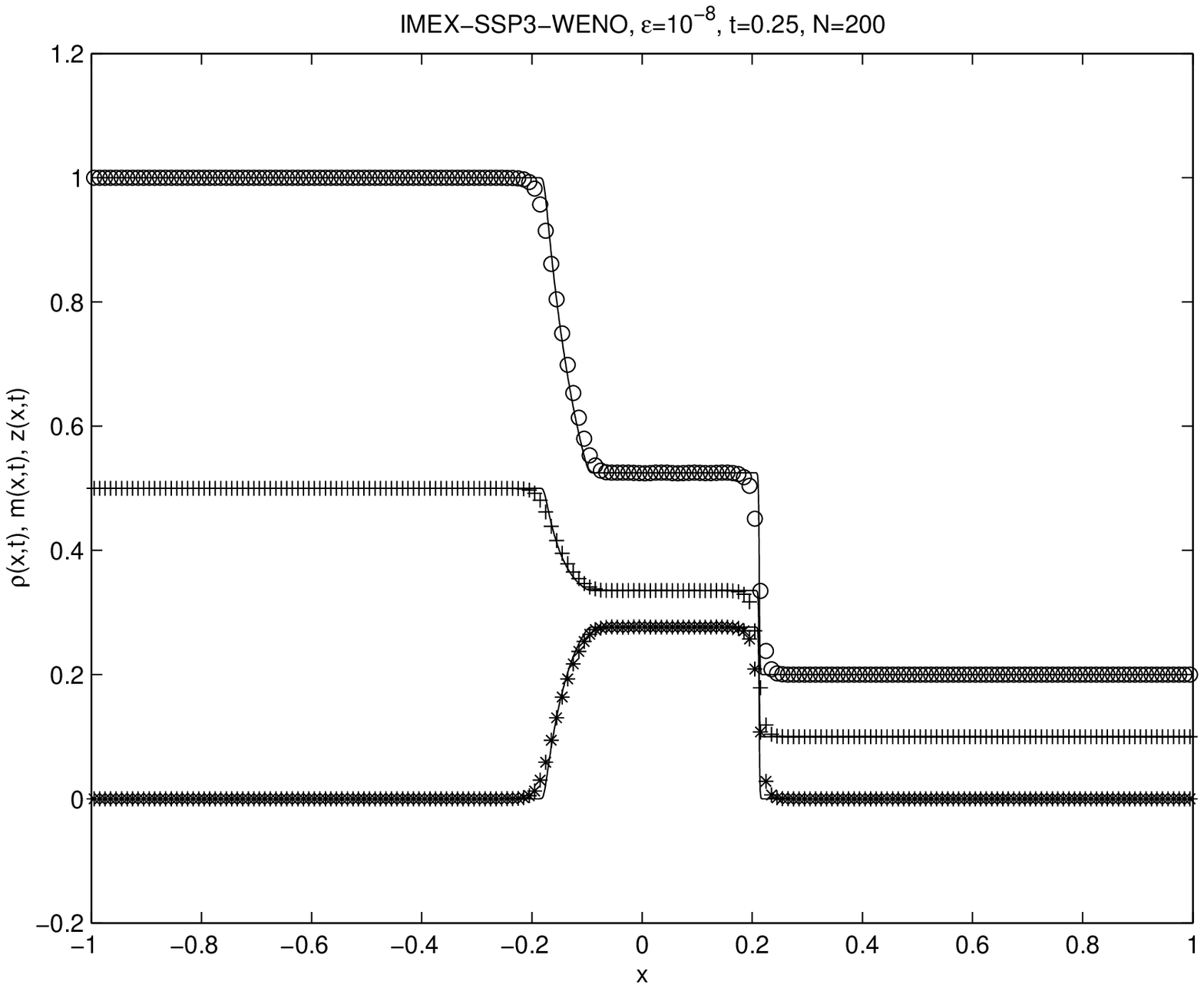}
} \vskip .5cm \caption{Numerical solution of the Broadwell
equations with initial data (\ref{eq:rim2}) for $\rho$($\circ$),
$m$($\ast$) and $z$($+$) at time $t=0.25$ for $\e=10^{-8}$. Left
IMEX-SSP2-WENO scheme, right IMEX-SSP3-WENO  scheme.}
\label{fig:2a}
\end{figure}

Next we test the shock capturing properties of the schemes in the
case of non smooth solutions characterized by the following two
Riemann problems \cite{CJR}
\bea
  \nonumber
  \rho_l=2,\qquad m_l=1,\qquad\qquad\quad z_l=1,\qquad x &<& 0.2, \\[-.3cm]
  \\
  \label{eq:rim1} \nonumber \rho_r=1,\qquad m_r=0.13962,\qquad
  z_r=1,\qquad x &>& 0.2, \eea \bea \nonumber
  \rho_l=1,\quad\qquad m_l=0,\qquad z_l=1,\qquad x &<& 0, \\[-.3cm]
  \\
  \label{eq:rim2} \nonumber \rho_r=0.2,\qquad m_r=0,\qquad
  z_r=1,\qquad x &>& 0.
\eea For brevity we report the numerical results obtained with the
second order IMEX-SSP2(2,2,2) and third order IMEX-SSP3(4,3,3)
schemes that we will refer to as IMEX-SSP2-WENO and IMEX-SSP3-WENO
respectively. The result are shown in Figures \ref{fig:1} and
\ref{fig:2a} for a Courant number $\Delta t/\Delta x = 0.5$. Both
schemes, as expected, give an accurate description of the solution
in all different regimes also using coarse meshes that do not
resolve the small scales. In particular the shock formation in the
fluid limit is well captured without spurious oscillations. We
refer to \cite{CJR, Jin-RK, LRR, Pa, RA} for a comparison of the
present results with previous ones.

\section{Applications}
Finally we present some numerical results obtained with
IMEX-SSP2-WENO and IMEX-SSP3-WENO concerning situations in which
hyperbolic systems with relaxation play a major role in
applications. The results have been obtained with $N=200$ grid
points. As usual the reference solution is computed on a much
finer grid.

\subsection{Shallow water}
First we consider a simple model of shallow water flow
\cite{Jin-RK}
\begin{eqnarray}
\nonumber
\partial_t h + \partial_x (h v)&=&0,\\[-.25cm]
\\
\nonumber
\partial_t (h v) + \partial_x (h+\frac12 h^2)&=&\frac{h}{\epsilon}(\frac{h}{2}- v),
\end{eqnarray}
where $h$ is the water height with respect to the bottom and $hv$
the flux.

The zero relaxation limit of this model is given by the inviscid Burgers
equation.

\begin{figure}[htb]
\centering{
\includegraphics[scale=0.4]{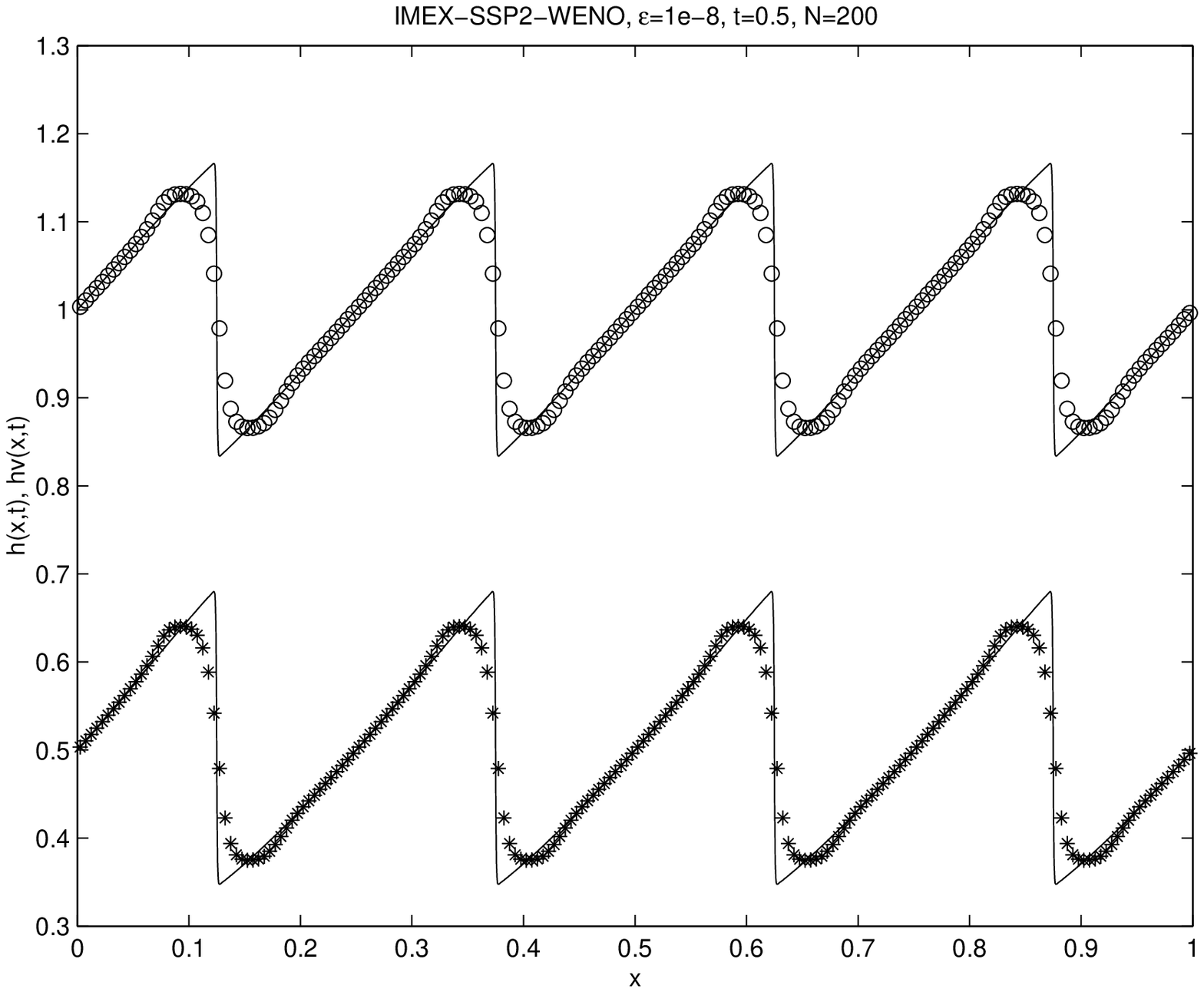}
\hskip .5cm
\includegraphics[scale=0.4]{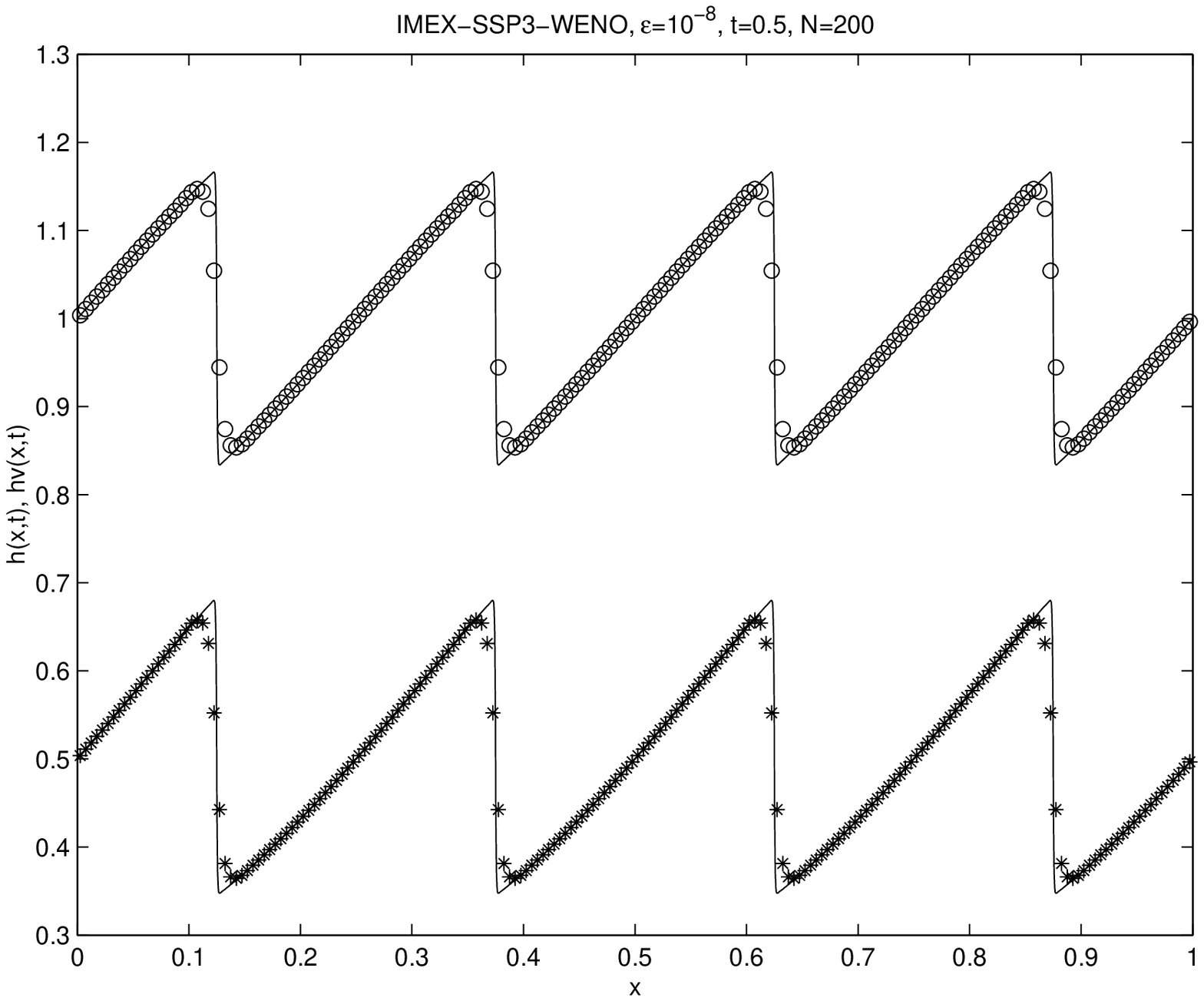}
} \vskip .5cm \caption{Numerical solution of the shallow water
model with initial data (\ref{eq:sh}) for $h$($\circ$) and
$hv$($\ast$) at time $t=0.5$ for $\e=10^{-8}$. Left IMEX-SSP2-WENO
scheme, right IMEX-SSP3-WENO scheme.} \label{fig:3}
\end{figure}

The initial data we have considered is \cite{Jin-RK}
\be
\nonumber
h = 1+0.2\sin(8\pi x), \quad
    hv = \frac{h^2}{2},
    \label{eq:sh}
\ee with $x\in [0,1]$. The solution at $t=0.5$ in the stiff regime
$\epsilon=10^{-8}$ using periodic boundary conditions is given in
Figure \ref{fig:3}. For IMEX-SSP2-WENO the dissipative effect due to
the use of the Lax-Friedrichs flux is very pronounced. As expected this
effect becomes less relevant with the increase of the order of
accuracy. We refer to \cite{Jin-RK} for a comparison with the present
results.

\subsection{Traffic flows}
In \cite{AR} a new macroscopic model of vehicular traffic has been presented.
The model consists of a continuity equation for the density $\rho$ of vehicles
together with an additional velocity equation that describes the mass flux
variations due to the road conditions in front of the driver. The model can be
written in conservative form as follows
\begin{eqnarray}
\nonumber
\partial_t \rho + \partial_x (\rho v)&=&0,\\[-.25cm]
\\
\nonumber
\partial_t (\rho w) + \partial_x (v\rho w)&=&A\frac{\rho}{T}(V(\rho)-v),
\end{eqnarray}
where $w=v+P(\rho)$ with $P(\rho)$ a given function describing the
anticipation of road conditions in front of the drivers and
$V(\rho)$ describes the dependence of the velocity with respect to
the density for an equilibrium situation. The parameter $T$ is the
relaxation time and $A>0$ is a positive constant.

If the relaxation time goes to zero, under the subcharacteristic
condition
\[
-P'(\rho) \leq V'(\rho) \leq 0, \quad \rho > 0,
\]
we obtain the Lighthill-Whitham \cite{Whitham} model
\begin{equation}
\partial_t \rho + \partial_x (\rho V(\rho)) = 0.
\end{equation}
A typical choice for the function $P(\rho)$ is given by
\[
P(\rho)= \left\{
\begin{array}{cc}
  \frac{c_v}{\gamma}\left(\frac{\rho}{\rho_m}\right)^\gamma & \gamma > 0,\\
  {c_v}\ln\left(\frac{\rho}{\rho_m}\right) & \gamma =
  0,
\end{array}
\right.
\]
where $\rho_m$ is a given maximal density and $c_v$ a constant
with dimension of velocity.
%This choice leads to the following
%natural choice for the equilibrium velocity $V(\rho)$
%\[
%V(\rho)=-c(P(\rho)-P(\rho_m)), \quad 0\leq c \leq 1.
%\]

In our numerical results we assume $A=1$ and an equilibrium
velocity $V(\rho)$ fitting to experimental data \cite{AKMR}
\[
V(\rho)=v_m
\frac{\pi/2+\arctan\left(\alpha\frac{\rho/\rho_m-\beta}{\rho/\rho_m-1}\right)}{\pi/2+\arctan\left(\alpha\beta\right)}
\]
with $\alpha=11$, $\beta=0.22$ and $v_m$ a maximal speed. We
consider $\gamma=0$ and, in order to fulfill the subcharacteristic
condition,  assume $c_v=2$. All quantities are normalized so that
$v_m=1$ and $\rho_m=1$.

\begin{figure}[htb]
\centering{
\includegraphics[scale=0.4]{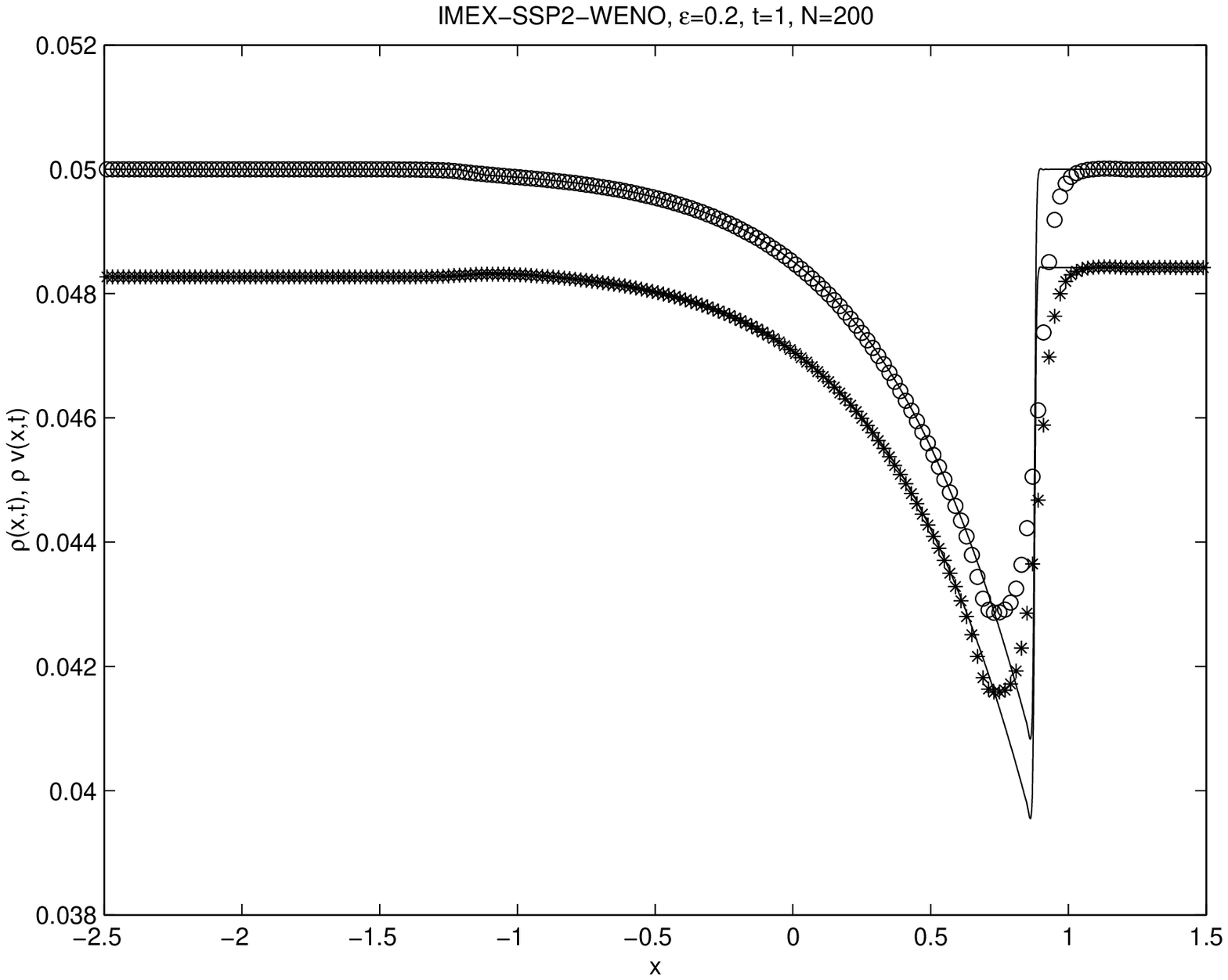}
\hskip .5cm
\includegraphics[scale=0.4]{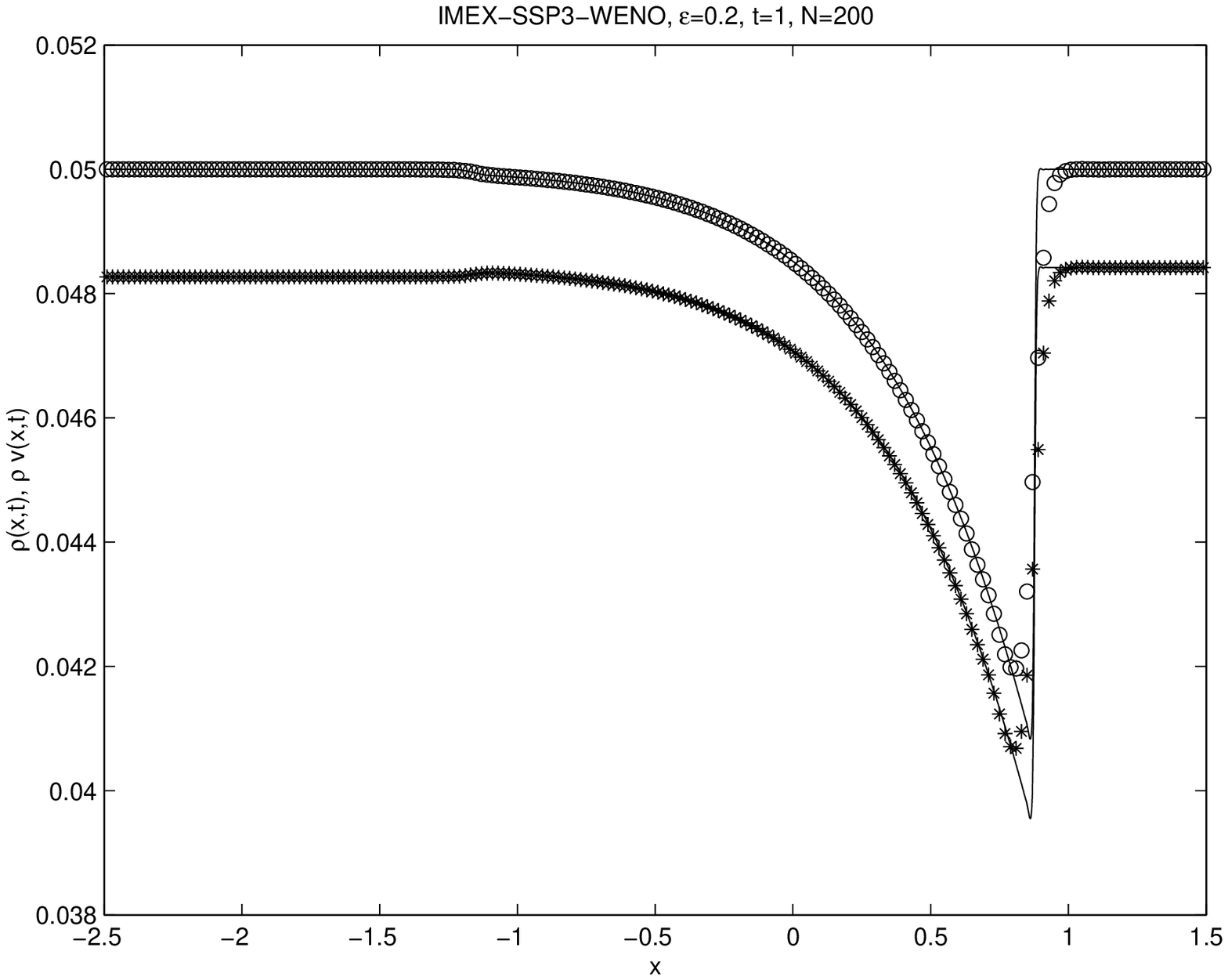}
} \vskip .5cm \caption{Numerical solution of the traffic model
with initial data (\ref{eq:tr1}) for $\rho$($\circ$) and $\rho
v$($\ast$) at time $t=1$ for $\e=0.2$. Left IMEX-SSP2-WENO scheme,
right IMEX-SSP3-WENO scheme.} \label{fig:tr1}
\end{figure}

We consider a Riemann problem centered at $x=0$ with left and
right states \be \rho_L=0.05,\quad v_L =0.05, \qquad
\rho_R=0.05,\quad v_R=0.5. \label{eq:tr1}\ee The solution at $t=1$
for $T=0.2$ is given in Figure \ref{fig:tr1}. The figure shows the
development of the density of the vehicles. Both schemes gives
very similar results. Again, in the second order scheme the shock
is smeared out if compared to the third order case. See \cite{AKMR}
for more numerical results.

\subsection{Granular gases}

We consider the continuum equations of Euler type for a granular
gas \cite{JR, To}. These equations have ben derived for a dense
gas composed of inelastic hard spheres. The model reads
\begin{eqnarray}
\nonumber
\rho_t + (\rho u)_x &=& 0,\\
\nonumber
(\rho u)_t + (\rho u^2 + p)_x &=& \rho g,\\[-.25cm]
\label{eq:NS}
\\[-.25cm]
\nonumber \left(\frac12\rho u^2 + \frac32\rho T \right)_t +
\left(\frac12 \rho u^3+ \frac32 u\rho T + pu\right)_x &=&
-\frac{(1-e^2)}{\epsilon}G(\rho){\rho^2 T^{3/2}},
\end{eqnarray}
where $e$ is the coefficient of restitution, $g$ the acceleration
due to gravity, $\epsilon$ a relaxation time, $p$ is the pressure
given by
\[
p=\rho T (1+2(1+e)G(\rho)),
\]
and $G(\rho)$ is the statistical correlation function. In our
experiments we assume
\[
G(\rho)=\nu\left(1-\left(\frac{\nu}{\nu_M}\right)^{\frac43\nu_M}\right)^{-1},
\]
where $\nu=\sigma^3\rho\pi/6$ is the volume fraction, $\sigma$ is
the diameter of a particle, and $\nu_M=0.64994$ is 3D random
close-packed constant.

We consider the following initial data \cite{CM} on the interval
$[0,10]$ \be \rho = 34.37746770,\quad
    v = 18, \quad P = 1589.2685472, \label{eq:gg1}\ee
which corresponds to a supersonic flow at Mach number $M_a=7$ (the
ratio of the mean fluid speed to the speed of sound). Zero-flux
boundary condition have been used on the bottom (right) boundary
whereas on the top (left) we have an ingoing flow characterized by
(\ref{eq:gg1}).

\begin{figure}[htb]
\centering{
\includegraphics[scale=0.41]{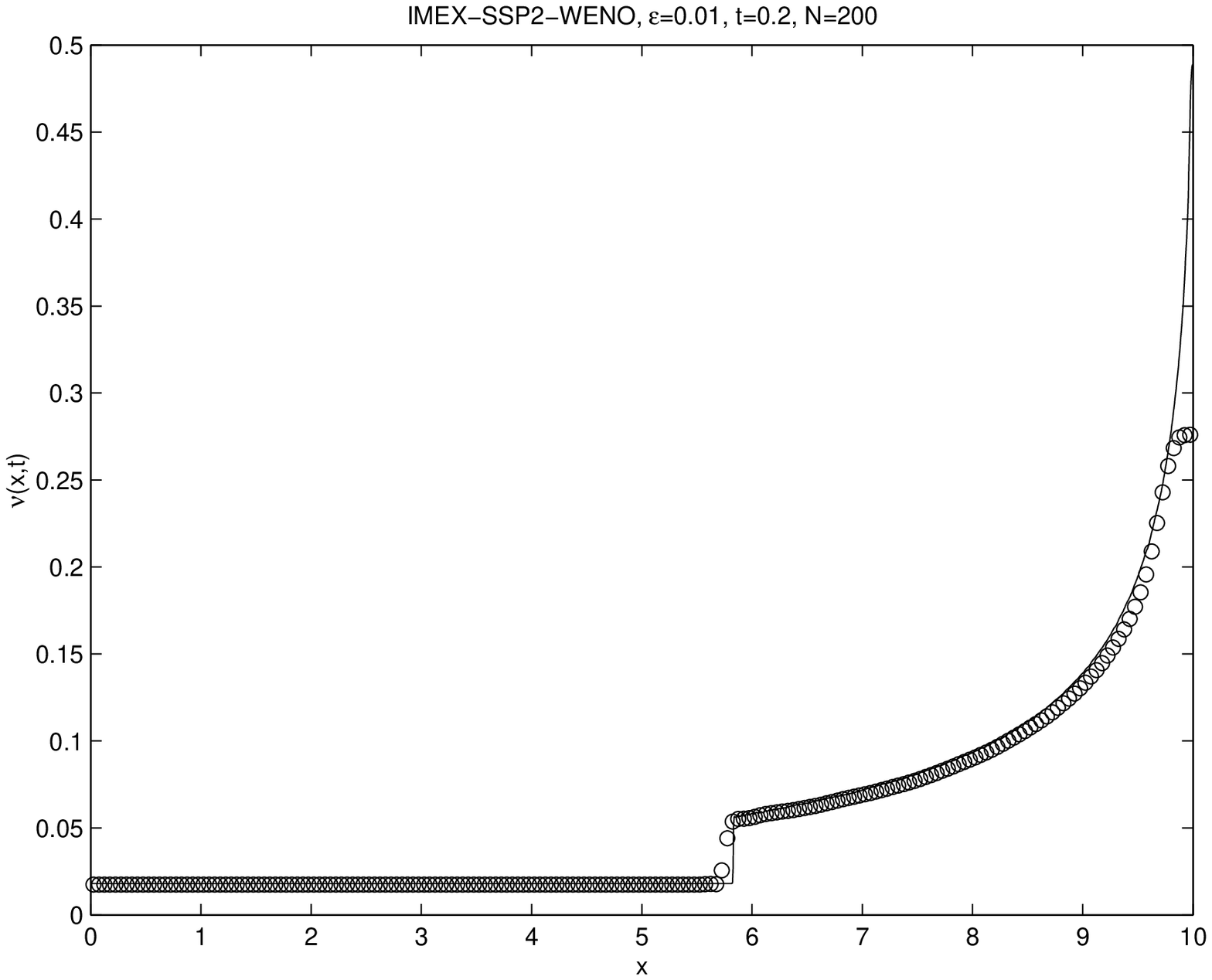}
\hskip .5cm
\includegraphics[scale=0.41]{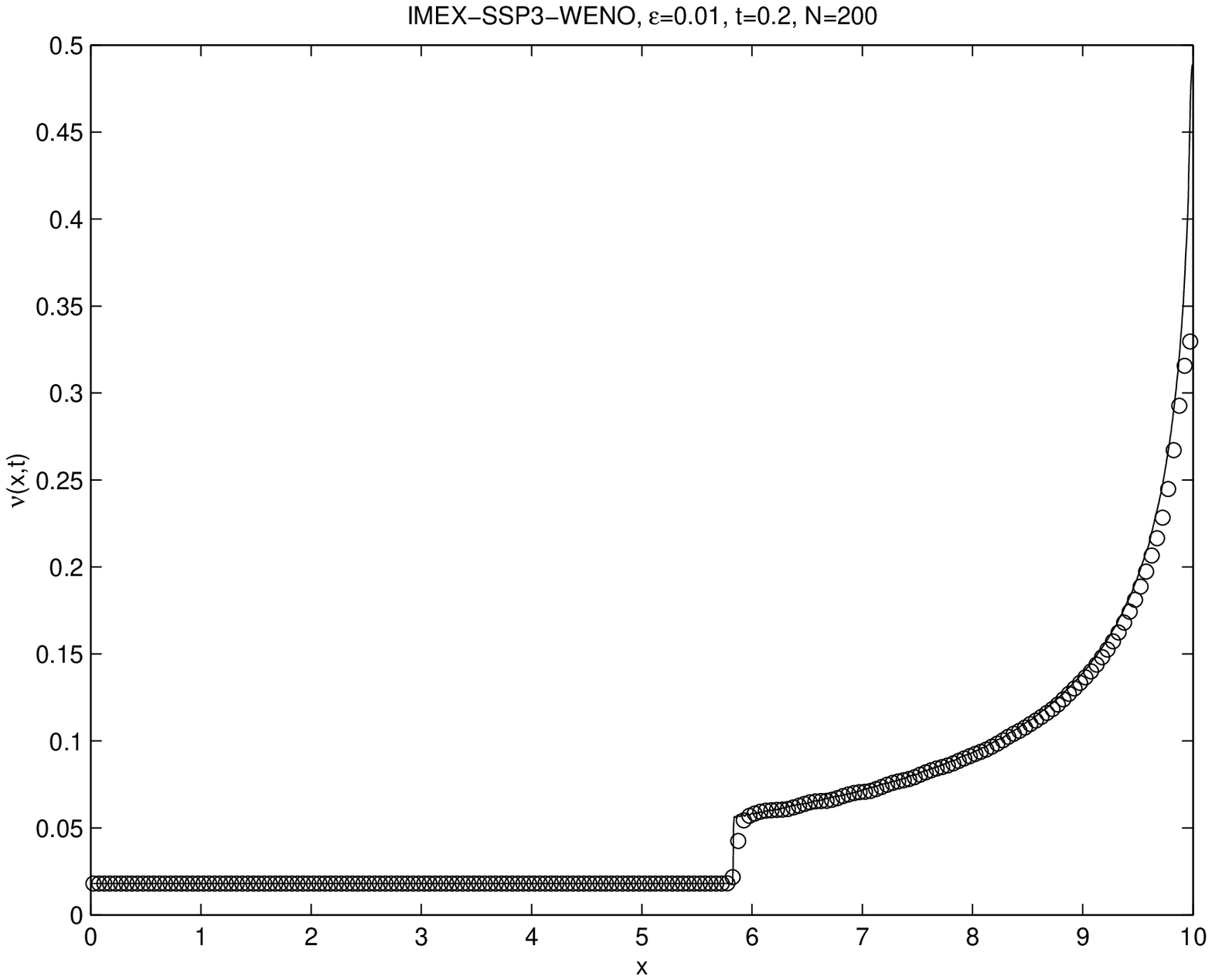}
\includegraphics[scale=0.41]{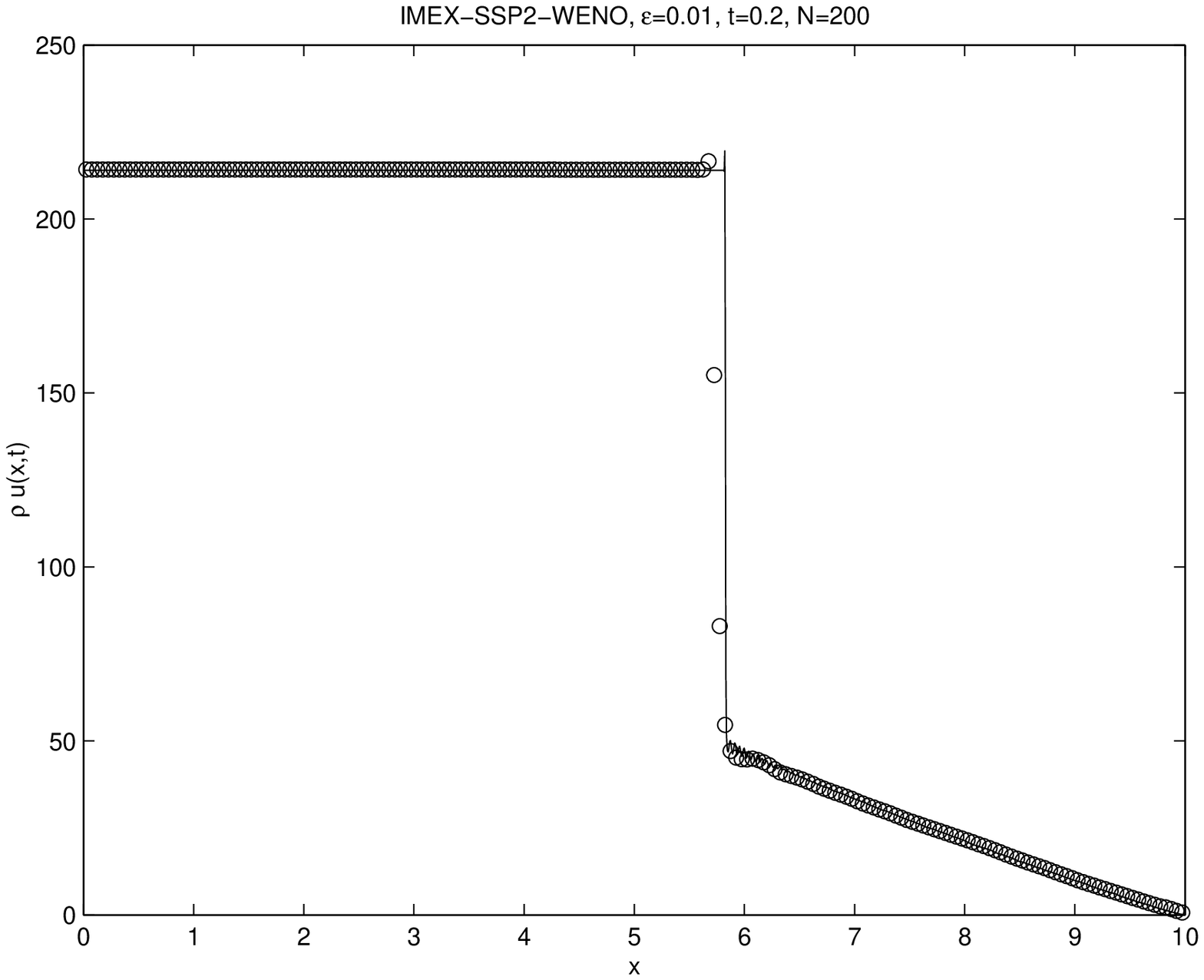}
\hskip .5cm
\includegraphics[scale=0.41]{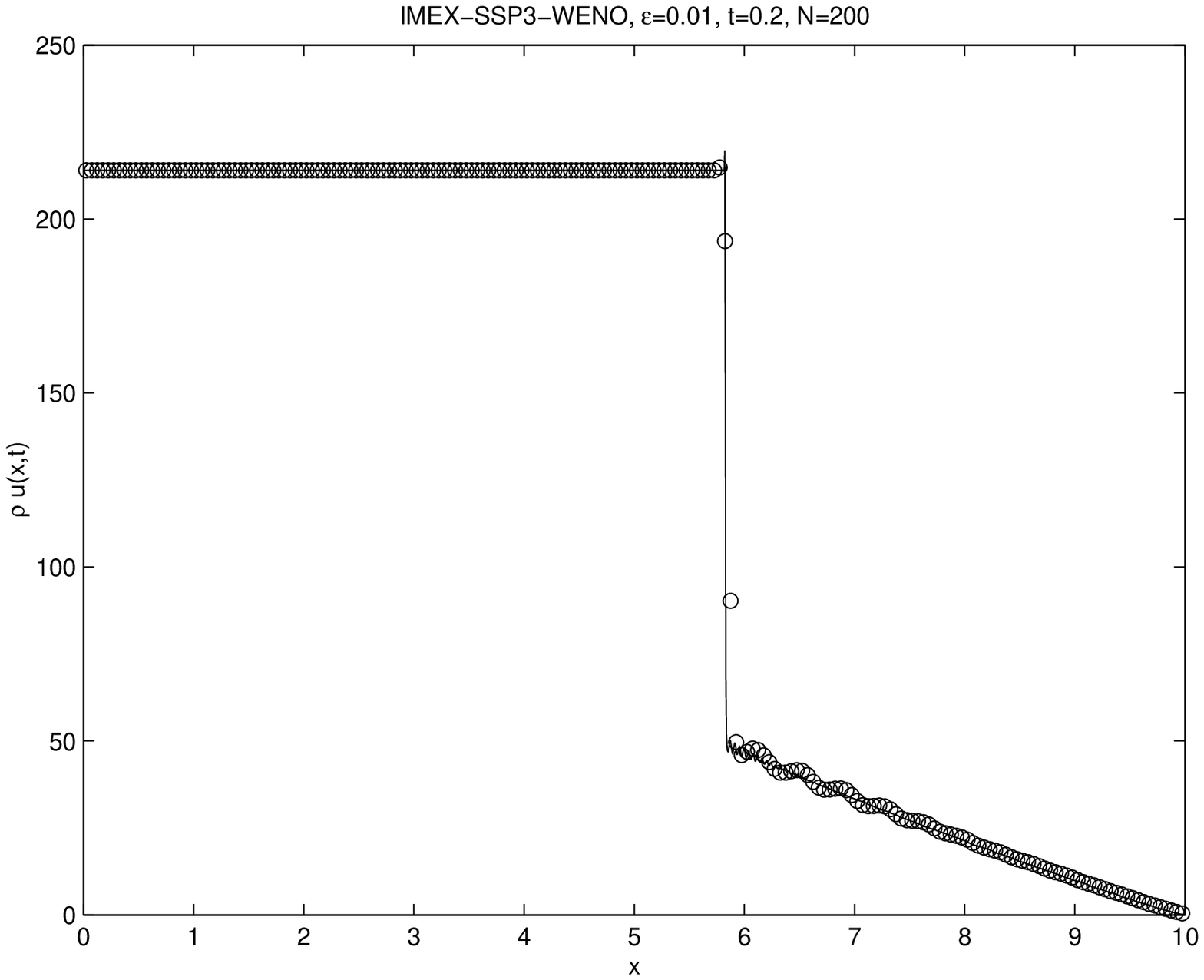}
\includegraphics[scale=0.41]{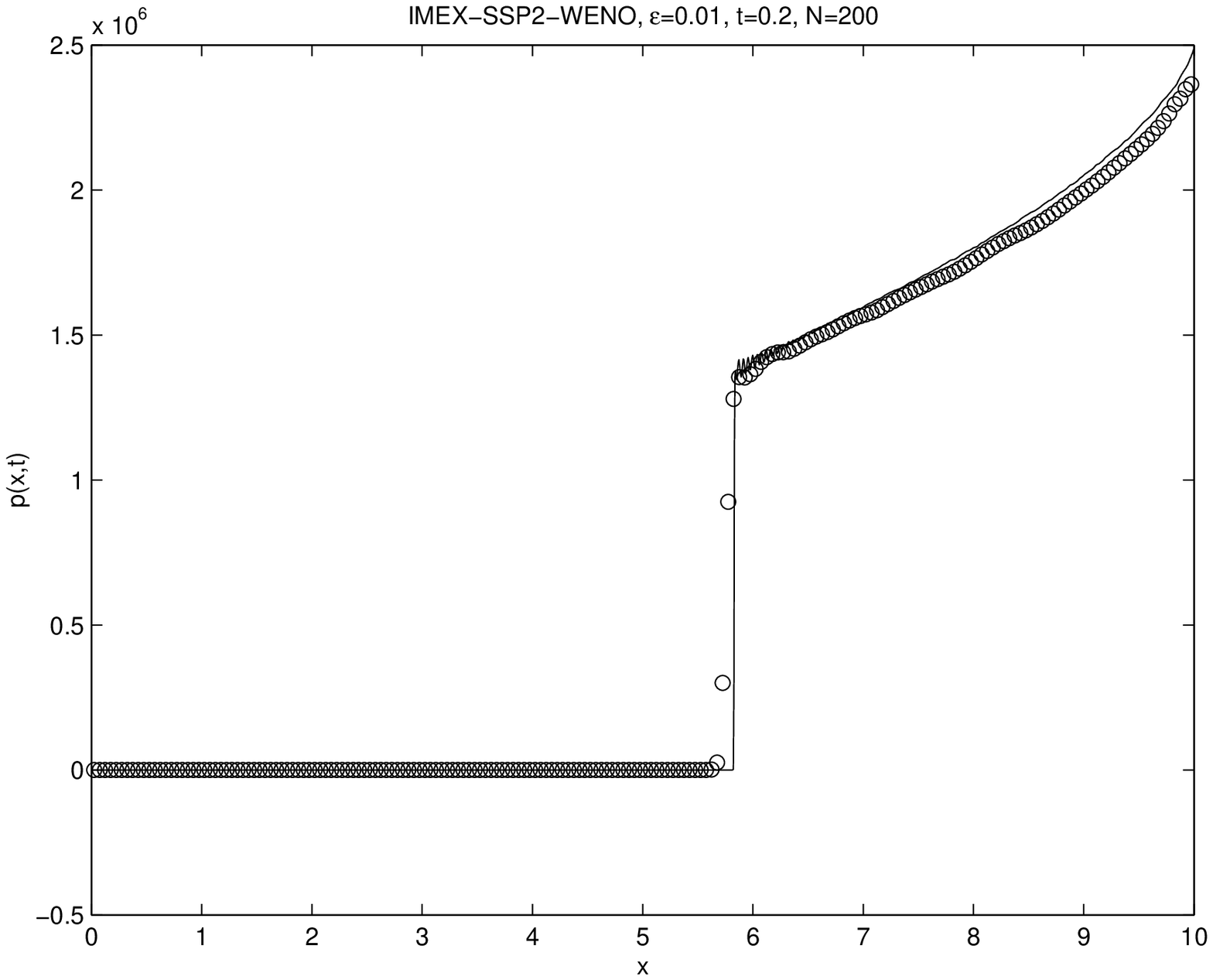}
\hskip .5cm
\includegraphics[scale=0.41]{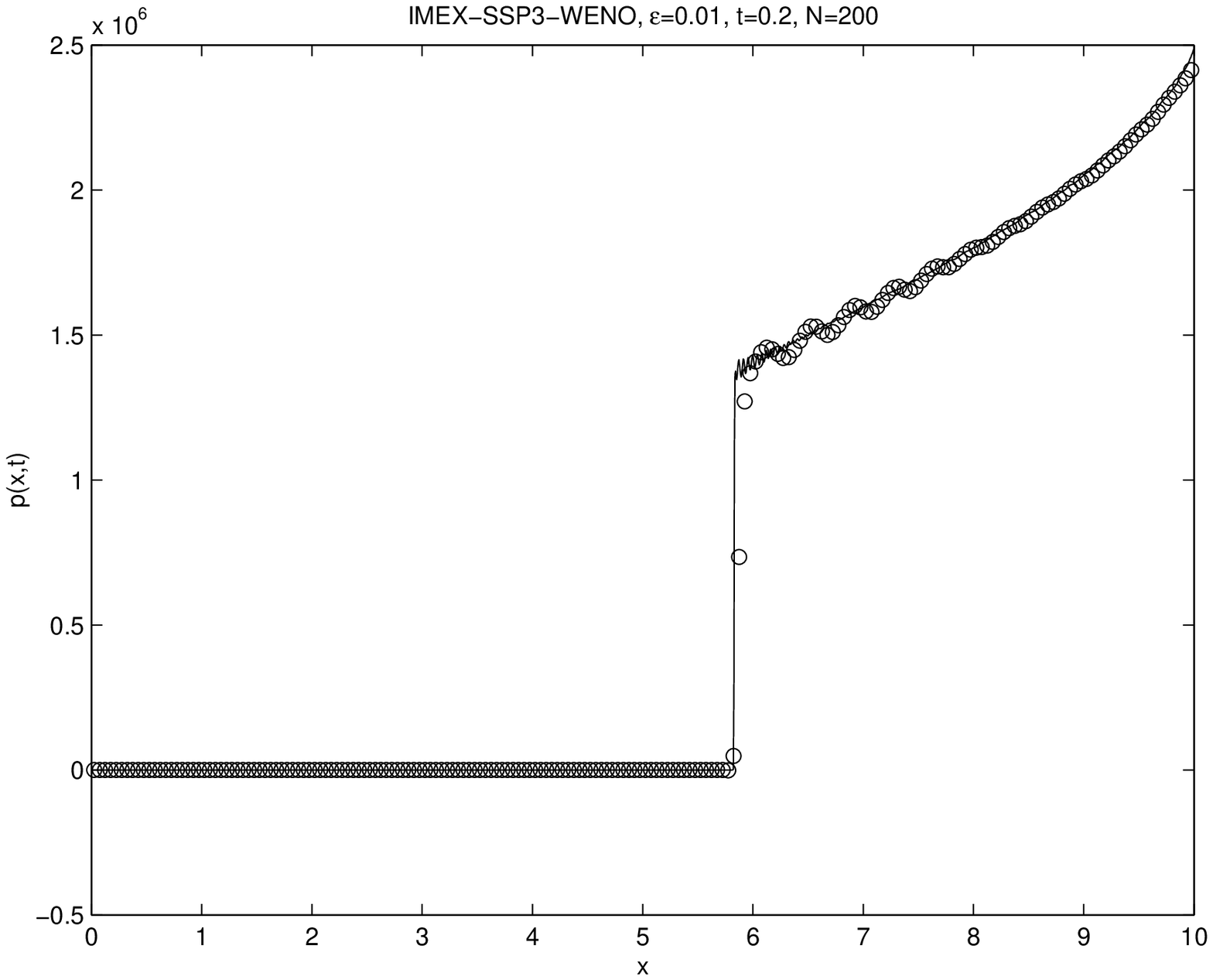}
} \vskip .5cm \caption{Numerical solution of the hydrodynamical
model of a granular gas with initial data (\ref{eq:gg1}). Left
column IMEX-SSP2-WENO scheme, right column IMEX-SSP3-WENO scheme.
From top to bottom, mass fraction $\nu$, velocity $\rho u$ and
pressure $p$.}\label{fig:g1}
\end{figure}

The values of the restitution coefficient and the particle
diameter have been taken $e=0.97$ and $\sigma=0.1$. We report the
solution at $t=0.2$ with $\epsilon=0.01$ in Figure \ref{fig:g1}
(see \cite{CM} for similar results). Due to the nonlinearity of
the source term the implicit solver has been efficiently solved
using Newton's method in each cell thanks to the use of finite
difference space discretization.

Both methods provide a good description of the shock that
propagates backward after the particles impact with the bottom.
Note that the second order method provides excessive smearing of
the layer at the right boundary. This problem is not present in
the third order scheme. However due to the use of conservative
variables we can observe the presence of small spurious
oscillations in the pressure profile.

\section*{Acknowledgements}
The authors would like to thanks the unknown referees for their
useful remarks and for pointing out reference \cite{Lubich}.


\begin{thebibliography}{99}
\bibitem{ACM} {G.~Akridis, M.~Crouzeix, C.~Makridakis}, {\em Implicit-explicit
multistep methods for quasilinear parabolic equations}, Numer.
Math. 82 (1999), no. 4, pp.~521--541.

\bibitem{RA}  M.~Arora, P.L.~Roe, {\it Issues and strategies for hyperbolic problems
with stiff source terms} in Barriers and challenges in computational fluid dynamics,
Hampton, VA, 1996, Kluwer Acad. Publ., Dordrecht, (1998), pp.~139--154.

\bibitem{AP98} U.~Ascher,, L.~Petzold, {\it Computer Methods for Ordinary Differential
Equations, and Differential Algebraic Equations}, SIAM,
Philadelphia, (1998).

\bibitem{ARS} {U.~Ascher, S.~Ruuth, R.~J.~Spiteri},  {\em Implicit-explicit
   Runge-Kutta methods for time dependent Partial Differential Equations\/},
   Appl. Numer. Math. 25, (1997), pp.~151--167.

\bibitem{ARW} {U.~Ascher, S.~Ruuth, B.~Wetton}, {\em Implicit-explicit
methods for time dependent PDE's\/}, SIAM J. Numer. Anal., 32,
(1995), pp.~797--823.

\bibitem{AR} {A.~Aw, M.~Rascle}, {\em Resurrection of second order models of traffic flow
?}, SIAM. J. Appl. Math. 60, (2000), pp.~916--938.

\bibitem{AKMR} {A.~Aw, A.~Klar, T.~Materne, M.~Rascle}, {\em Derivation of continuum traffic flow models from microscopic
follow the leader models}, SIAM J. Appl. Math. 63, (2002),
pp.~259--278.

\bibitem{Butcher} {J.~Butcher}, {\em The numerical analysis of Ordinary differential
   equations. Runge-Kutta and general linear methods.\/}. John Wiley \& Sons, Chichester and
   New York (1987).

\bibitem{NBR} {M.~Briani, R.~Natalini, G.~Russo}, {\em
Implicit--Explicit numerical schemes for jump diffusion
processes}, in preparation (2004).

\bibitem{CJR} {R.~E.~Caflisch, S.~Jin, G.~Russo}, {\em Uniformly accurate schemes
   for hyperbolic systems with relaxation\/}, SIAM J. Numer. Anal.,
   34, (1997), pp.~246--281.

\bibitem{CK} {C.~A.~Kennedy, M.~H.~Carpenter},
{\em Additive Runge-Kutta schemes for convection-diffusion-reaction equations},
 Appl. Numer. Math. 44 (2003), 139--181.

\bibitem{Cecchi-Redivo-Russo} {M.~Cecchi Morandi, M.~Redivo-Zaglia, G.~Russo},
   {\em Extrapolation methods for hyperbolic systems with relaxation\/},
   Journal of Computational and Applied Mathematics 66 (1996), pp.~359--375.

\bibitem{CLL} {G.~Q.~Chen, D.~Levermore, T.~P.~Liu}, {\it Hyperbolic conservations laws
with stiff relaxation terms and entropy}, Comm. Pure Appl. Math.,
47, (1994), pp.~787--830.

\bibitem{Dekker-Verwer} {K.~Dekker, J.~G.~Verwer}, {\em Stability of Runge-Kutta Methods
   for Stiff Nonlinear Differential Equations\/}, North-Holland, Amsterdam (1984).

\bibitem{Dia-Shatzman1} {B.~O.~Dia, M.~Schatzman},
   {\em Estimation sur la formule de Strang\/}, C. R. Acad. Sci. Paris, t.320,
   S\'erie I (1995), pp.~775--779.

\bibitem{Dia-Shatzman2} {B.~O.~Dia, M.~Schatzman},
   {\em Commutateur de certains semi-groupes holomorphes et applications
   aux directions altern\'ees\/}, Mathematical Modelling and Numerical
   Analysis 30 (1996), pp.~343--383.

\bibitem{FHV} {J.~Frank, W.~H.~Hundsdorfer, J.~G.~Verwer}, {\em On the stability of
implicit-explicit linear multistep methods}, Special issue on time integration
(Amsterdam, 1996). Appl. Numer. Math. 25 (1997), no. 2-3, 193--205.

\bibitem{GS} {S.~Gottlieb, C.~-W.~Shu},
   {\em Total Variation Diminishing Runge-Kutta
   schemes\/}, Math. Comp. 67 (1998), pp.~73--85.

\bibitem{GST} {S.~Gottlieb, C.~-W.~Shu, E.~Tadmor}, {\em Strong-stability-preserving high order time discretization
methods\/}, SIAM Review, 43 (2001), pp.~89--112.

\bibitem{Hairer} {E.~Hairer}, {\em Order conditions for numerical methods for Partitioned ordinary
   differential equations\/}, Numerische Mathematik 36 (1981) pp.~431-445.

\bibitem{HNW} {E.~Hairer, S.~P.~N{\o}rsett, G.~Wanner},
   {\em Solving ordinary differential equations, Vol.1 Nonstiff
    problems}, Springer-Verlag,     New York (1987).

\bibitem{HW} {E.~Hairer, G.~Wanner},
   {\em Solving ordinary differential equations, Vol.2 Stiff and
    differential-algebraic problems}, Springer-Verlag,     New York (1987).

\bibitem{Hofer} {E.~Hofer} {\em A partially implicit method for large stiff systems
   of Ode's with only few equations introducing small time-constants\/}, SIAM J. Numer. Anal.
   13, (1976) pp.~645-663.

\bibitem{Lubich} {T.~Jahnke, C.~Lubich}, {\it Error bounds for exponential operator
splitting}, BIT, (2000), pp.~735--744.

\bibitem{JR} {J.~Jenkins, M.~Richman},
{\em Grad's 13-moment system for a dense gas of inelastic spheres},
Arch. Rat. Mechanics, 87, (1985), pp.~355--377.

\bibitem{Jin-RK} {S.~Jin}, {\em Runge-Kutta methods for hyperbolic systems
   with stiff relaxation terms\/} J. Comput. Phys., 122 (1995), pp.~51--67.

\bibitem{JinLe} {S.~Jin, C.~D.~Levermore},
{\em Numerical Schemes for hyperbolic conservation laws with stiff
relaxation terms}, J. Comp. Physics, 126 (1996), pp.~449--467.

\bibitem{JX} S.~Jin, Z.P.~Xin, {\em The relaxation schemes for systems
of conservation laws in arbitrary space dimensions}, Comm. Pure
Appl. Math. 48 (1995), no. 3, 235--276.

%\bibitem{JPT} {  S.~Jin, L.~Pareschi and G.~Toscani} {\em Diffusive relaxation schemes for multiscale discrete-velocity kinetic equations}, SIAM J. Num. Anal., 35 (1999) pp.~2405-2439.

\bibitem{leveque:numerical-methods} {R.~J.~LeVeque},
   {\em Numerical Methods for Conservation Laws},
   Lectures in Mathematics, Birkhauser Verlag, Basel (1992).

\bibitem{LRR} {S.~F.~Liotta, V.~Romano, G.~Russo},
   {\em Central schemes for balance laws of relaxation type\/}, SIAM J. Numer. Anal., 38, (2000), pp.~1337--1356.

\bibitem{Liu} {T.~P.~Liu}, {\it Hyperbolic conservation laws with relaxation}, Comm. Math. Phys., 108, (1987), pp.~153--175.

\bibitem{CM} A.~Marquina, S.~Serna, {\em Capturing shock waves in inelastic granular gases},
UCLA-CAM Report, 04-04 february (2004).

\bibitem{Minion} {M.~L.~Minion}, {\em Semi-implicit spectral deferred correction
methods for ordinary differential equations}, Comm. Math.
Sciences, 1, No.3, p.~471--500.

%\bibitem{NP} {G.~Naldi, L.~Pareschi},
%{\em Numerical schemes for hyperbolic systems with stiff diffusive relaxation}, SIAM J. Numer. Anal. to appear.

\bibitem{MR}  I.~M\"uller, T.~Ruggeri, {Rational extended thermodynamics}, Springer-Verlag, Berlin, (1998).

\bibitem{Pa} {L.~Pareschi}, {\em Central differencing based numerical schemes for hyperbolic conservation
laws with stiff relaxation terms}, SIAM J. Num. Anal., 39, (2001),
pp.~1395-1417.

\bibitem{PR} {L.~Pareschi, G.~Russo}, {\em Implicit-explicit {R}unge-{K}utta
schemes for stiff systems of
  differential equations}, {Advances Theo. Comp. Math.}, 3, (2000), pp.~269--289.

\bibitem{PR2} {L.~Pareschi, G.~Russo}, {\em High order asymptotically
strong-stability-preserving methods for hyperbolic systems with
stiff relaxation}, {Proceedings HYP2002}, Pasadena USA, Springer
(2003), pp.~241--255.

\bibitem{PR3} {L.~Pareschi, G.~Russo}, {\em Stability analysis of
Implicit-Explicit Runge-Kutta schemes}, preprint (2004).

\bibitem{Shu-central} Jianxian Qiu, Chi-Wang Shu, {\em On the construction, comparison, and local
characteristic decomposition for high-order central WENO schemes}, J. Comput.
Phys. 183 (2002), no. 1, 187--209

\bibitem{Shu-TVD} {C.~-W.~Shu}, {\em Total variation diminishing time discretizations},
SIAM J. Sci. Stat. Comput., 9, (1988), 1073--1084.

\bibitem{Shu-CIME} {C.~-W.~Shu}, {\em Essentially Non Oscillatory and Weighted
    Essentially NOn OScillatory Schemes for Hyperbolic Conservation Laws}, in
  Advanced numerical approximation of nonlinear hyperbolic equations, Lecture
  Notes in Mathematics, 1697, (2000).

\bibitem{Shu-OsherII} Chi-Wang Shu, S. Osher,
{\em Efficient implementation of essentially nonoscillatory
shock-capturing schemes}, J. Comput. Phys. 77, (1988), no. 2,
pp.~439--471.

\bibitem{Sod} G.A.~Sod, {\it A survey of several finite difference methods
for systems of nonlinear hyperbolic conservation laws}, J. Comp.
Phys. 27, (1978), pp.~1--31.

\bibitem{SR} {R.~J.~Spiteri, S.~J.~Ruuth}, {\em A new class of optimal
strong-stability-preserving time discretization
methods}, SIAM. J. Num. Anal. 40, (2002), no. 2, 469--491.

\bibitem{Strang} {G.~Strang}, {\em On the
construction and comparison of difference schemes\/}, SIAM J.
Numer. Anal. {5}, (1968) pp.~505--517.

\bibitem{tadmor:approximate} {E.~Tadmor}, {\em Approximate Solutions of Nonlinear
   Conservation Laws}, Cockburn, Johnson, Shu and Tadmor Eds.,
   Lecture Notes in Mathematics, N.~1697 (1998).

\bibitem{To}
G. Toscani,
\newblock Kinetic and hydrodinamic models of nearly elastic granular flows,
\newblock {\em Monatsch. Math.} (to appear)

\bibitem{Whitham} {G.~B.~Whitham}, {\it Linear and nonlinear waves}, Wiley, New York, (1974).

\bibitem{Zhong} {X.~Zhong}, {\em Additive Semi-Implicit Runge-Kutta methods for computing high speed nonequilibrium reactive
flows}, J. Comp. Phys. 128,(1996),  pp.~19--31.

\end{thebibliography}
\end{document}